\documentclass[preprint,12pt]{elsarticle}

\usepackage{amsmath,amssymb,amsthm}
\usepackage{latexsym}            

\usepackage{epsfig}
\usepackage{epstopdf}
\usepackage{amsfonts,amscd}
\usepackage[ruled,vlined]{algorithm2e}
\usepackage{color,comment}
\usepackage{amsbsy,bm}
\usepackage{subfigure}
\usepackage[utf8]{inputenc}          
\usepackage[T1]{fontenc}

\newcommand\vel{\mathbf{u}}
\newcommand\bv{\mathbf{v}}
\newcommand\bx{\mathbf{x}}

\newcommand\bn{\mathbf{n}}

\newcommand\bD{\boldsymbol{D}}

\newcommand\btau{\boldsymbol{\tau}}
\newcommand\bsigma{\boldsymbol{\sigma}}

\newcommand\bj{\boldsymbol{j}}

\newcommand{\be}{\begin{eqnarray}}
\newcommand{\ee}{\end{eqnarray}}
\newcommand{\ben}{\begin{eqnarray*}}
	\newcommand{\een}{\end{eqnarray*}}

\newcommand {\na} {{\nabla}}
\newtheorem{thm}{Theorem}[section]

\newtheorem{lma}[thm]{Lemma}

\newtheorem{rmk}{Remark}[section]
\allowdisplaybreaks \allowdisplaybreaks[4]

\begin{document}
	
	\begin{frontmatter}
		
		
		
		\title{An energy stable $C^0$ finite element scheme for a quasi-incompressible phase-field model of moving contact line with variable density}
		
		\author[mymainaddress]{Lingyue Shen}
		\author[mythirdaddress,myfourthaddress]{Huaxiong Huang}
		\author[mymainaddress]{Ping Lin}	
		\author[mythirdaddress,myfifthaddress]{Zilong Song}
		\author[myfourthaddress,mysecondaryaddress]{Shixin Xu  \corref{mycorrespondingauthor}}
		\cortext[mycorrespondingauthor]{Corresponding author}
		\ead{shixin.xu@dukekunshan.edu.cn}
		\address[mymainaddress]{Department of Mathematics, University of Dundee, Dundee DD1 4HN, Scotland, UnitedKingdom.}
		\address[mysecondaryaddress]{Duke Kunshan University, 8 Duke Ave, Kunshan, Jiangsu, China. }
		\address[mythirdaddress]{Department of Mathematics and Statistics, York University, Toronto, Ontario, Canada. }
		\address[myfourthaddress]{Centre for Quantitative Analysis and Modelling, Fields Institute for Research in Mathematical Sciences, Toronto, Ontario, Canada}
		\address[myfifthaddress]{Department of Mathematics, University of California, Riverside, 900 University Avenue, Riverside, CA, USA}
		\begin{abstract}
			In this paper, we focus on modeling and simulation of two-phase flow problems with  moving contact lines and variable density. A thermodynamically consistent   phase‐field model with General Navier Boundary Condition is developed based on the concept of quasi-incompressibility and the energy variational method.  A mass conserving   C0 finite element scheme is proposed to solve the PDE system.  Energy stability is achieved at the fully discrete level. Various numerical results confirm that the proposed scheme  for both P1 element and P2 element are energy stable. 
		\end{abstract}
		
		
		
		
	\end{frontmatter}
	
	\section{Introduction}
	The modeling and simulation of moving contact lines (MCLs), where the interface  of two or more immersible fluids intersects with a solid wall \cite{dussan_v._motion_1974,dussan_spreading_1979},  have attracted much attention in recent years.  Applications of MCLs in industries and medical fields (for example, printing \cite{kumar_liquid_2015}, spray cooling of surfaces \cite{breitenbach_drop_2018}, blood clot \cite{xu_model_2017}, microfluidics \cite{tung_contact_2013}, surfactant \cite{xu2014level,zhang_derivation_2014}) have motivated scientific interests and mathematical challenges on associated issues such as the stress singularity and contact angle hysteresis.  In order to  model the dynamics around the contact lines, various types of models and approximations have been developed, such as direct molecular dynamics simulations \cite{koplik_molecular_1988, koplik_molecular_1989,smith_moving_2018}, phase-filed models \cite{bao2012finite,jiang_numerical_2015,qian_variational_2006,qian_molecular_2006,salgado_diffuse_2013,shen_efficient_2015},  microscopic–macroscopic hybrid model \cite{hadjiconstantinou_hybrid_1999,ren_heterogeneous_2005}, front tracking model \cite{lai_numerical_2010,ren2007boundary,ren2010continuum,zhang_derivation_2014} and Lattice Boltzmann model \cite{huang2018lattice}. 
	For reviews of the current status of the MCLs problem, we refer to the articles  \cite{bonn_wetting_2009} and \cite{snoeijer_moving_2013}. 
	
	Among those models, phase-field method (or diffusive interface method) \cite{anderson_diffuse-interface_1998,gurtin_two-phase_1996,jacqmin_calculation_1999} is one of the most popular and powerful methodologies. It has two main advantages.  Firstly, it is easy to track the interface  and numerically implement even if there are topological changes \cite{du_retrieving_2005}. Secondly, it can be  derived by energy-based variational approach \cite{feng2005energetic,yue_diffuse-interface_2004,liu_decoupled_2015}. As a result, the obtained system is compatible with the law of energy dissipation, which makes it possible to design efficient and energetically stable numerical schemes.
	
	One of the main challenges in phase-field method is to model the immersible two-phase flow with different densities. When the density ratio between the two phases is small, it could be  handled by the Boussinesq approximation \cite{liu2003phase}. However, it could not be extended to the case with a large density ratio due to its underlying assumption \cite{liu_decoupled_2015}.  One key problem arises from the inconsistency between the mass conservation and the incompressibility especially near the diffusive interface region.
	It was first pointed out by Lowengrub \cite{lowengrub_quasiincompressible_1998} and later by Shen et al. \cite{liu_decoupled_2015,shen_mass_2013}.   Two  main approaches are proposed to overcome this difficulty: one is based on volume  averaged velocity; the other is based on the mass averaged velocity. For the volume averaged velocity model, the incompressiblity is assumed everywhere including the interfacial region \cite{abels_thermodynamically_2010,boyer_nonhomogeneous_2001,ding_diffuse_2007,shen_energy_2010,shen_phase-field_2010}. An thermodynamically consistent and frame invariant model was developed by Abels et al. \cite{abels_thermodynamically_2012}, where the mass conservation equation is modified with a mass correction term. On the other hand, for the mass averaged velocity method,  the mass conservation is assured instead of incompressibility. This naturally yields the quasi-incompressible Navier-Stoker-Cahn-Hilliard (q-NSCH) model \cite{lowengrub_quasiincompressible_1998,guo_mass_2017}, which in fact leads to a slightly compressible mixture only inside the interfacial region.   
	
	In the present paper, we first rederive and generalize the thermodynamcially consistent q-NSCH model in \cite{lowengrub_quasiincompressible_1998} from a variational point of view by combining with  the Energy Variational Approach (EnVarA)  \cite{yue_diffuse-interface_2004,eisenberg2010energy,liu2019energetic} and Onsager's Variation Principle   \cite{zhang_derivation_2014,ren2007boundary,ren2010continuum,abels_thermodynamically_2012}.  It  starts from two functionals for the total energy and dissipation, together with the kinematic equations based on physical laws of conservation.  The specific forms of the fluxes and stresses in the kinematic equations could be obtained by taking the time derivative of the total energetic functional and comparing with the predefined dissipation functional. More details could be found in \cite{xu2018osmosis}. In addition to bulk energy and dissipation, the energy and dissipation on the boundary are introduced to model the dynamics of contact lines. Our energy variational approach consistently yields both the correct bulk equations (the q-NSCH system)  and  a modified General Navier-Stokes Boundary Condition (GNBC) for the case of mass averaged velocity. The density effect on the contact line is explicitly modeled compared with the traditional  GNBC \cite{qian_variational_2006,qian_molecular_2006,gao_gradient_2012,gao_efficient_2014,xu2018three,yu_numerical_2017} in the case of volume averaged velocity, where the effect is modeled implicitly by the bulk and boundary interactions.


	The second goal of our paper  is to design an efficient energy stable scheme for the obtained q-NSCH system with large density ratio. There are not many such schemes developed for the MCLs. For the incompressible NSCH system,  the development of such schemes may be found in
	\cite{salgado_diffuse_2013,gao_gradient_2012,gao_efficient_2014,luo2017efficient,zhang2016phase, 
		dong_time-stepping_2012,shen_decoupled_2015,yang2018efficient,yu_numerical_2017} and only a few of them \cite{gao_efficient_2014,yu_numerical_2017,dong_time-stepping_2012} are for variable density MCL models using the volume-averaged velocity (satisfying the incompressible condition in the whole domain). We shall develop an energy stable scheme for our thermodynamically consistent variable density q-NSCH system using the mass-averaged velocity.
	Based on the author's  previous works \cite{jiang_numerical_2015,guo_mass_2017,guo2014numerical}, we design a mass conservative  $C^0$ finite element method for the q-NSCH system with a consistent discrete energy law.  Thanks to a $\Delta p$ term in the quasi-incompressible condition, which is similar to the  pressure stabilization of pseudo-compressibility methods \cite{brezzi1984stabilization,rannacher1993numerical,shen1992pressure}, q-NSCH system 
	does not need to satisfy the Babuska-Brezzi
	inf-sup condition \cite{shen1992pressure,lin1997sequential,lin2004finite,lin2006simulations}. This may be considered as another benefit of our  quasi-incompressible NSCH system.
	
	
	The rest of paper is organized as follows. In Section \ref{section:model}, we present the thermodynamically consistent derivation of the q-NSCH system and its non-dimensionlization. The $C^0$ finite element algorithm for the q-NSCH system and the energy stable analysis are shown in Section \ref{section:numerial}. Section \ref{section:simulation} presents the numerical results, including the convergence case study,  and the examples of moving droplets and  rising bubbles.  

	\section{Mathematical Model}\label{section:model}
	\subsection{Mass-averaged velocity and laws of conservation}
\begin{figure}[!ht]
		\centering
		\subfigure[] {
			\scalebox{0.4}[0.4]{\includegraphics{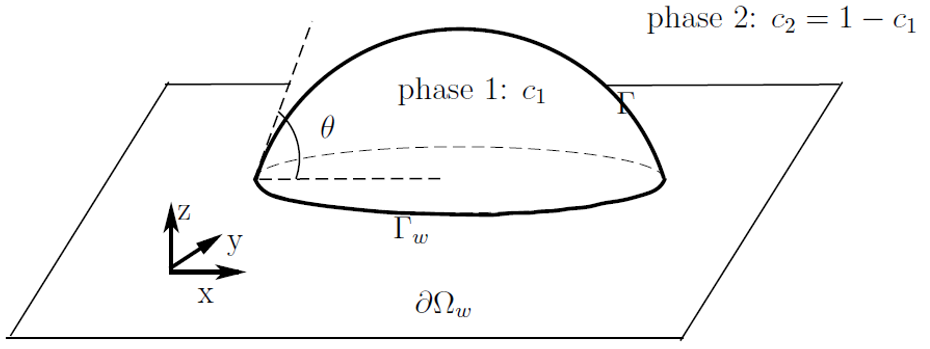}}
		}
		\subfigure[] {
			\scalebox{0.4}[0.4]{\includegraphics{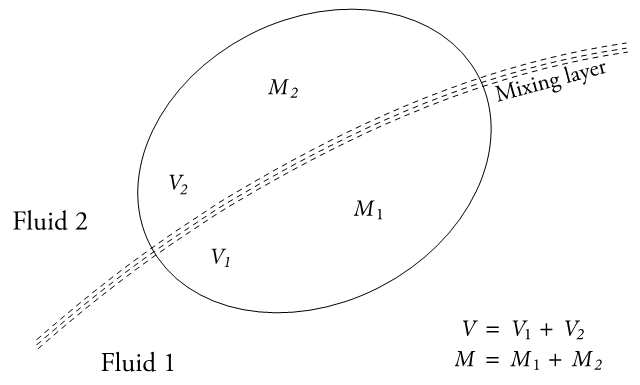}}
		}
		\caption{ Schematic of  moving contact line problems (a) and interface (b). 
			\label{Schematic} }
	\end{figure}

	We consider a complex mixture consisting of two phase fluids with different densities. The interface of two fluids intersects with the wall $\partial\Omega_w$ at the contact line $\Gamma_w$ (see Fig. \ref{Schematic} (a)). Around the interface, we choose a control volume $V(t)$, where there are two phases labeled by $i=1,2$ with volume $V_i$ and mass $M_i$ (see Fig. \ref{Schematic} (b)).
	If the  local average density of each phase  is denoted by $\bar{\rho}_i=M_i/V$ and pure phase density   is denoted by $\rho_i=M_i/V_i$, then  density of mixture is
	\be\label{rhov}
	\rho=\frac{M}{V}=\frac{M_1}{V}+\frac{M_2}{V}=\bar{\rho}_1+\bar{\rho}_2.
	\ee
	Let  $c_i=\frac{M_i}{M}$ be the mass fraction of each phase. Then we have \cite{guo2015thermodynamically}
	\be\label{rhoc}
	\frac{1}{\rho}=\frac{V}{M}=\frac{V_1}{M}+\frac{V_2}{M} = \frac{c_1}{\rho_1}+\frac{c_2}{\rho_2}=\frac{c} {\rho_1}+\frac{1-c}{\rho_2},
	\ee
	where $c=c_1$ is adopted in the last equality.
	\begin{rmk}
		Note that according to the definition \eqref{rhoc}, the mixture density $\rho$ is almost constant  everywhere except  in the interfacial region. 
	\end{rmk}

	If we assume those two fluids move  with velocities $\vel_i$ ($i=1,2$), then the mass conservation of each phase inside the control volume is 
	\begin{equation}
	\frac{\partial\bar{\rho}_i}{\partial t} +\nabla\cdot(\vel_i\bar{\rho}_i) = 0, \quad i=1,2.
	\end{equation}
	Introducing the mass averaged velocity as 
	\begin{equation}
	\rho\vel = \bar\rho_1\vel_1+\bar\rho_2\vel_2,
	\end{equation}
	and combining with Eq.\eqref{rhov} yields the conservation of mass for the mixture 
	\begin{equation}\label{con_mass}
	\frac{\partial\rho}{\partial t} +\nabla\cdot(\vel\rho) = 0.  
	\end{equation}

	Next, with an arbitrary volume $V(t)\in\Omega$, laws of conservation state
	\begin{eqnarray}
	\frac{d}{dt}\int_{V(t)}\rho c dx = -\int_{\partial V(t)}\bj_c \cdot\bn dS\label{con_c},\\
	\frac{d}{dt}\int_{V(t)}\rho \vel dx = \int_{\partial V(t)}(\bsigma_{\eta}+\bsigma_{c})\cdot\bn dS\label{con_v}.
	\end{eqnarray}
	Here, the first equation is the conservation of phase-field function (phase 1) and  $\bj_{ c}$ is the flux of phase-field function.  The second equation is the conservation of momentum where $\bsigma_{\eta}$ is the viscous stress and $\bsigma_{ c}$ is the extra stress induced by  two-phase  interface due to nonzero $\nabla c$.
	
	Thanks to the Reynolds transport theory \cite{guo_mass_2017,xu2018osmosis}, Eqs.\eqref{con_mass}-\eqref{con_v}  yield the following kinematic equations in the domain $\Omega$
	\be\left\{\begin{array}{l}\label{assumption}
		\rho\frac{Dc}{D t}  = -\nabla\cdot \bj_{ c},\\
		\frac{D\rho}{dt}+\rho\nabla\cdot\vel=\frac{\partial\rho}{\partial t} +\nabla\cdot(\vel\rho) = 0,\\
		\rho\frac{D\vel}{D t}  = \nabla\cdot\bsigma_{\eta} +\nabla\cdot\bsigma_{ c},
	\end{array}\right.\ee
	where $\frac{D}{Dt} = \frac{\partial }{\partial t}+\vel\cdot\nabla$ is the material derivative.
	
	By definition of $\rho = \rho( c(\mathbf{x},t))$  in Eq. \eqref{rhoc},    above equations yield the quasi-incompressibility condition \cite{lowengrub_quasiincompressible_1998,guo_mass_2017}
	\be\label{quasi_incom}
	\nabla\cdot \vel=-\frac{1}{\rho}\frac{D\rho}{Dt} = \frac{1}{\rho^2} \frac{d  \rho}{d c}  \left( \nabla\cdot \bj_{ c}  \right).
	\ee
	In the present case, we denote
	\be\label{quasi_incom2}
	\alpha = -\frac{1}{\rho^2} \frac{d  \rho}{d c} = \frac{\rho_2-\rho_1}{\rho_1\rho_2},
	\ee
	then the quasi-incompressibility condition is written as
	\be\label{quasi_incom3}
	\nabla\cdot \vel =-\alpha \nabla\cdot \bj_{ c}.
	\ee
	\begin{rmk}
		Equations \eqref{quasi_incom2}-\eqref{quasi_incom3} show that the quasi-incompressibilty condition depends on the density difference of two fluids. When the two phases have the same density, i.e. $\rho_1=\rho_2$, it will  consistently degenerate to the incompressibility condition. It makes a difference when two fluids have large density ratio and near the interfacial region \cite{lowengrub_quasiincompressible_1998,guo_mass_2017,guo2015thermodynamically}. 
		
	\end{rmk}
	On the boundary of domain $\partial\Omega$, the following boundary conditions are used
	\be\left\{\begin{array}{l}\label{assumption_bd}
		\vel\cdot\bn =0, ~~\vel^{s}\cdot \btau_i= u^s_{\tau_i}=f_{\tau_i}\\
		\frac{D_{\Gamma}c}{Dt}= J_{\Gamma},\\
		\bj_{ c}\cdot \bn = 0,\\
	\end{array}\right.\ee
	where $\vel^{s} = \vel_{\tau}-\vel_w$ with $\vel_{\tau} = \vel-(\vel\cdot\bn)\bn$ is the fluid slip velocity with respect to the wall,  $\frac{D_{\Gamma}c}{Dt}=\frac{\partial  c}{\partial t} + \vel\cdot\nabla_{\Gamma} c$ is the surface material derivative, the Allen-Cahn type boundary condition is used for $ c$ and $\nabla_{\Gamma} = \nabla-\bn(\bn\cdot\nabla)$ is surface gradient on the boundary $\partial\Omega$. The quantities $f_{\tau_i}$ and $J_{\Gamma}$ are to be determined.  During the derivation, we  assume the solid wall is fixed, i.e. $\vel^{s} = \vel_{\tau}$.

	\subsection{Model derivation} 
	
	Now we start to derive the exact forms of $\bj_{ c}$, $\bsigma_{\eta}$, $\bsigma_{ c}$ in Eq. (\ref{assumption}), $f_{\tau_i}$ and $J_{\Gamma}$ in Eq. (\ref{assumption_bd}) by using energy variational method. 
	
	The total energy consists of the kinetic energy, the phase mixing energy and the energy on solid wall boundary $\partial\Omega_w$
	\be\label{energy}
	E^{tot} &= &E_{kin}+E_{mix}+E_{w}\nonumber\\
	&=& \int_{\Omega}\frac{\rho(c)|\vel|^2}{2} d\bx +\int_{\Omega}\lambda_{ c}\rho(c)\left(G( c) +\frac{\gamma^2}{2}|\nabla c|^2\right)d\bx
	+\int_{\partial\Omega_w}f_w( c)dS,
	\ee
	together with
	\be
	G( c)=\frac 1 4  c^2(1- c)^2, \quad f_w( c) = -\frac{\sigma}{2}\cos(\theta_s) \sin(\frac{(2 c-1)\pi}{2}),
	\ee
	where $\lambda_{ c}$ is the  mixing energy density, $\gamma$ is  the capillary width  of the interface, $\theta_s$  is static contact angle and $\sigma$ is surface tension. The mixing energy $E_{mix}$ represents the competition between a homogeneous bulk mixing energy density term $G( c)$ (`hydrophobic' part) that enforces total separation of the two phases into pure components, and a gradient distortional term $\frac{|\na c|^2}2 $  (`hydrophilic' part) that represents the nonlocal interactions between two phases and penalizes spatial heterogeneity.
	

	The dissipation functional is composed of the dissipation due to fluid friction and irreversible mixing of two phases in bulk and the dissipation on the boundary
	\be\label{dissipation}
	\Delta = \int_{\Omega}2\eta( c)|\bD_{\eta}|^2d\bx+\int_{\Omega}\lambda( c)|\nabla\cdot\vel|^2d\bx  +\int_{\Omega}\frac 1 {\mathcal{M}}|\bj_{c}|^2d\bx+\int_{\partial\Omega_w}\left( \frac{1}{\mathcal{M}_{\Gamma}}J_\Gamma^2 +\beta_\Gamma |\vel^s|^2\right)dS,
	\ee
	where $\lambda(c)$ and $\eta(c)$ are the two Lam\'e coefficients, $\bD_{\eta}= (\nabla\vel+(\nabla\vel)^T)/2$ is the strain rate, $\mathcal{M}$ is mobility coefficient in bulk, $\mathcal{M}_{\Gamma}$ is mobility coefficient on the wall, $\beta_{\Gamma}(c)$ is wall friction coefficient. 
	In the present paper,  $\eta(c)$  and $\beta_{\Gamma}$ are  approximated by 
	$$
	\frac1 {\eta(c)} =\frac c {\eta_1} + \frac{(1-c)}{\eta_2},
	~~\frac1 {\beta_{\Gamma}(c)} =\frac {c} {\beta_{\Gamma_1}} + \frac{(1-c)}{\beta_{\Gamma_2}},
	$$
	where   $\eta_i$ and $\beta_{\Gamma,i}$ with $i=1,2$ are coefficients of each phase.  
	

	During the derivation, the following  lemma is frequently used.
	\begin{lma}\label{transportlma}
		For a continuous function $f(\bx,t)$, if the density $\rho$ satisfies the conservation law \eqref{con_mass} in the domain $\Omega$ and $\vel\cdot \bn=0$ on the boundary $\partial\Omega$, then we have 
		\be
		\frac{d}{dt}\int_{\Omega}\rho(\bx,t) f(\bx,t)d\bx = \int_{\Omega}\rho\frac{Df}{Dt}d\bx.\nonumber
		\ee
	\end{lma}

	By taking the  time derivative of the total energetic functional, we have
	\be\label{dEdt}
	\frac{dE^{tot}}{dt} = \frac{d}{dt}E_{kin}+\frac{d}{dt}E_{mix} +\frac{d}{dt}E_{w}= I_1+I_2+I_3.
	\ee
	For the first term in (\ref{dEdt}), using the last two equations in Eq.\eqref{assumption} yields
	\be\label{I_1}
	I_1 &=&\frac{d}{dt}\int_{\Omega}\frac{\rho|\vel|^2}{2}d\bx \nonumber\\
	&=&-\int_{\Omega}(\bsigma_{\eta}:\nabla\vel +\bsigma_{ c}:\nabla\vel)d\bx+\int_{\Omega}\alpha\nabla  p  \cdot \bj_{ c}d\bx-\int_{\Omega} p  \nabla\cdot \vel d\bx\nonumber\\
	&&+\int_{\partial\Omega_w}((\bsigma_{\eta}+\bsigma_{ c})\cdot\bn)\cdot\vel_{\tau} dS,
	\ee
	where we have introduced a Lagrangian multiplier $p$ with respect to the quasi-compressibility condition (\ref{quasi_incom}) and have used the boundary conditions $\vel\cdot\bn = 0$ and $\bj_{ c}\cdot\bn=0$. For the second term in (\ref{dEdt}), using the first equation in Eq.\eqref{assumption} and last two boundary conditions in Eq. \eqref{assumption_bd} yields
	\be\label{I_2}
	I_2& =&\frac{d}{d t}\int_{\Omega}\rho\lambda_{ c}\left(G( c)+\frac{\gamma^2}{2}|\nabla c|^2\right)d\bx\nonumber\\
	&=&\int_{\Omega}\nabla\mu\cdot\bj_{ c} d\bx-\int_{\Omega}\lambda_{ c} \gamma^2(\rho\nabla c \otimes\nabla c):\nabla\vel d\bx+\int_{\partial\Omega_w}\rho \lambda_{ c}\gamma^2\partial_n c \frac{D_{\Gamma} c}{Dt}dS,
	\ee
	where $\mu = \lambda_{ c} \left(\frac{dG}{d c}-\frac{1}{\rho}\gamma^2\nabla\cdot(\rho\nabla c)\right)$.
	The detailed derivations of Eqs.\eqref{I_1} -\eqref{I_2} are given in Appendix \ref{energyvariation}. The last term $I_3$ in (\ref{dEdt}) yields
	\be\label{I_3}
	I_3 = \frac{d}{dt}\int_{\partial\Omega_w}f_wdS = \int_{\partial\Omega_w}\frac{d f_w}{d c}\frac{\partial c}{\partial t}dS.
	\ee
	
	Combining Eqs.\eqref{I_1} -\eqref{I_3}, we obtain the derivative of the energy functional 
	\be\label{eq20}
	\frac{dE^{tot}}{dt} &=& -\int_{\Omega}\bsigma_{\eta}:\nabla\vel d\bx-\int_{\Omega}\left(\bsigma_{ c}+ \lambda_{ c}\gamma^2\rho\nabla c\otimes\nabla c\right):\nabla\vel d\bx\nonumber\\
	&&+\int_{\Omega}\nabla\mu\cdot\bj_ c d\bx +\int_{\Omega}\nabla (\alpha  p)  \cdot \bj_{ c}d\bx-\int_{\Omega} p\nabla\cdot \vel d\bx\nonumber\\
	&&+\int_{\partial\Omega_w}((\bsigma_{\eta}+\bsigma_{ c})\cdot\bn)\cdot\vel_{\tau} dS+\int_{\partial\Omega_w}\rho \lambda_{ c}\gamma^2\partial_n c \frac{D_{\Gamma} c}{Dt}dS+\int_{\partial\Omega_w}\frac{d f_w}{d c}\frac{\partial c}{\partial t}dS\nonumber\\
	&=& -\int_{\Omega}\bsigma_{\eta}:\nabla\vel d\bx -\int_{\Omega}\left(\bsigma_{ c}+ \lambda_{ c}\gamma^2\rho\nabla c\otimes\nabla c\right):\nabla\vel d\bx\nonumber\\
	&&+\int_{\Omega}\nabla\tilde\mu\cdot\bj_ c d\bx-\int_{\Omega} p \nabla\cdot \vel d\bx\nonumber\\
	&&+\int_{\partial\Omega_w}\left((\bsigma_{\eta}+\bsigma_{ c})\cdot\bn-\frac{df_w}{d c}\nabla_{\Gamma} c\right)\cdot\vel_{\tau} dS+\int_{\partial\Omega_w}L( c) \frac{D_{\Gamma} c}{Dt}dS \nonumber\\
	&=& -\int_{\Omega}\bsigma_{\eta}:\nabla\vel d\bx -\int_{\Omega}\left(\bsigma_{ c}+ \lambda_{ c}\gamma^2\rho\nabla c\otimes\nabla c\right):\nabla\vel d\bx\nonumber\\
	&&+\int_{\Omega}\nabla\tilde\mu\cdot\bj_ c d\bx-\int_{\Omega} p \nabla\cdot \vel d\bx\nonumber\\
	&&+\int_{\partial\Omega_w}\left((\bsigma_{\eta}+\bsigma_{ c})\cdot\bn-\frac{df_w}{d c}\nabla_{\Gamma} c\right)\cdot\vel_{\tau} dS+\int_{\partial\Omega_w}L( c) J_\Gamma dS,
	\ee
	where we have defined
	\be\label{defofunkonwflux}
	\tilde\mu &=& \mu+\alpha p =\lambda_ c \frac{dG}{d c}-\frac {\lambda_{ c}\gamma^2} {\rho}\nabla\cdot(\rho\nabla c)+\alpha p ,\\
	L( c) & =& \rho\lambda_{ c}\gamma^2\partial_n c+\frac{df_w}{d c}.\label{eq22}
	\ee
	Using energy dissipation law   ${dE^{tot}}/{dt}= -\Delta$ \cite{eisenberg2010energy,xu2014energetic} and comparing (\ref{eq20}) with the predefined dissipation functional in Eq.\eqref{dissipation} yield 
	\be\label{eq23}
	\left\{\begin{array}{l}
		\bj_ c=-\mathcal{M}\nabla\tilde\mu,\\
		\bsigma_{\eta} = 2\eta\bD_{\eta}+\lambda\nabla\cdot\vel I -pI =\eta (\nabla\vel+(\nabla\vel)^T)+\lambda\nabla\cdot\vel I -pI ,\\
		\bsigma_{ c} =- \lambda_{ c}\gamma^2 \rho\left(\nabla c\otimes\nabla c\right),\\
		J_{\Gamma} = -\mathcal{M}_{\Gamma}L( c),\\
		u^s_{\tau_i} = \beta_{\Gamma}^{-1}\btau_i\cdot(-(\bsigma_{\eta}+\bsigma_{ c})\cdot\bn+\frac{df_w}{d c}\nabla_{\Gamma} c).
	\end{array}\right.
	\ee
	By the definition of $\bsigma_\eta$ and $\bsigma_{ c}$, the slip boundary condition (last equation in $(\ref{eq23})$) could be further written in the GNBC format
	\be
	u^s_{\tau_i} = \beta_{\Gamma}^{-1}\btau_i\cdot(-\bsigma_{\eta}\cdot\bn+L( c)\nabla_{\Gamma} c).
	\ee

	To summarize, we have the following model for the two-phase flow with variable density for three unknowns $ c,\vel, p$, in domain $\Omega$, 
	\begin{subequations}
		\label{model}
		\begin{eqnarray}
		&&\rho\frac{Dc}{D t}  =  \nabla\cdot(\mathcal{M}\nabla\tilde\mu),\label{phieq}\\
		&&\tilde\mu =\lambda_ c \frac{dG}{d c}-\frac {\lambda_{ c}\gamma^2} {\rho}\nabla\cdot(\rho\nabla c)+\alpha p,\label{mueq}\\
		&&\frac{D\rho}{Dt}+\rho\nabla\cdot\vel=0,\label{masseq}\\
		&&\rho\frac{D\vel}{D t}  = \nabla\cdot(2\eta \bD_{\eta}) + \nabla(\lambda\nabla\cdot\vel)-\nabla p -\nabla\cdot(\lambda_{ c}\gamma^2\rho\nabla c\otimes\nabla c),\label{nseq}
		\end{eqnarray}
	\end{subequations}
	with boundary conditions on $\partial\Omega$
	\be\left\{\begin{array}{l}\label{model_bd}
		\frac{D_{\Gamma} c}{Dt} = - \mathcal{M}_{\Gamma}L( c),\\
		\nabla\tilde\mu\cdot\bn = 0,\\
		\vel\cdot\bn = 0,\\
		u_{\tau_i}^s =\beta_{\Gamma}^{-1}(-(\bn\cdot\bsigma_{\eta}\cdot\btau_i)+L( c)\partial_{\tau_i} c), i=1,2.
	\end{array} 
	\right.
	\ee
	where $L(c)$ and $\bsigma_{\eta}$ are defined in \eqref{eq22} and \eqref{eq23}.
	\begin{rmk}
		Note that  in the above boundary conditions \eqref{model_bd}, the density effect on the contact line  dynamics is explicitly modeled both in the  boundary dynamics of phase-field and  in velocity slip boundary condition through $L(c)$ term.  
	\end{rmk}
	
	It is worth noting that the above system satisfies the following energy dissipation law.
	\begin{thm}\label{energylaw}
		If $ c, \vel, p$ are smooth solutions of above system \eqref{model}-\eqref{model_bd}, then
		the following energy law is satisfied:
		\be
		\frac{dE^{tot}}{dt} &\!\!\!=\!\!\!& \frac{d}{dt}\left\{ \int_{\Omega}\frac{\rho|\vel|^2}{2} d\bx +\int_{\Omega}\lambda_{ c}\rho\left(G( c) +\frac{\gamma^2}{2}|\nabla c|^2\right)d\bx
		+\int_{\partial\Omega_w}f_w( c)dS\right\}\nonumber\\
		&\!\!\!=\!\!\!&-\int_{\Omega}2\eta|\bD_{\eta}|^2d\bx-\int_{\Omega}\lambda|\nabla\cdot\vel|^2d\bx  -\int_{\Omega}\mathcal{M}|\nabla\tilde\mu|^2d\bx\nonumber\\
		&&-\int_{\partial\Omega_w}\left(\mathcal{M}_{\Gamma}\left|L( c)\right|^2 +\beta_{\Gamma}|\vel^{s}|^2\right)dS.
		\ee 
	\end{thm}
	\noindent \textbf{Proof:}  The main idea of the proof is obtained by multiplying the phase-field equation \eqref{phieq} by $\tilde{\mu}$, multiplying the chemical potential equation \eqref{mueq} by $\frac{d c}{dt}$, multiplying the mass conservation equation \eqref{masseq} by $p$, 
	multiplying the Navier-Stokes equation \eqref{nseq} by $\vel$, and summing them up.  
	
	Taking the inner product of the  phase-field equation \eqref{phieq} with $\tilde{\mu}$ results in the following equation
	\be
	\int_{\Omega}\rho\frac{D c}{Dt}\tilde{\mu}d\bx = -\int_{\Omega}\mathcal{M}|\nabla\tilde{\mu}|^2d\bx,
	\ee
	where we used the boundary condition $\partial_n\tilde{\mu}=0$ in \eqref{model_bd}.
	
	Multiplying the chemical potential \eqref{mueq} by $\rho\frac{D c}{Dt}$ yields
	\be
	\int_{\Omega}\tilde{\mu}\rho\frac{D c}{Dt} = \int_{\Omega}\lambda_{ c}\rho\frac{DG}{Dt}d\bx+\int_{\Omega}\lambda\gamma^2\nabla c\nabla\left(\frac{D c}{Dt}\right)d\bx-\int_{\Omega}\frac{p}{\rho}\frac{d \rho}{dc}\frac{D c}{Dt}d\bx-\int_{\partial\Omega_w}\lambda\gamma^2\rho\partial_n c \frac{Dc}{Dt}dS.
	\ee 
	Summing up the above two equations, we have 
	\be\label{pheenergy}
	&&\int_{\Omega}\lambda_{ c}\rho\frac{DG}{Dt}d\bx+\int_{\Omega}\lambda\gamma^2\nabla c\nabla\left(\frac{Dc}{Dt}\right)d\bx\nonumber\\
	&=&-\int_{\Omega}\mathcal{M}|\nabla\tilde{\mu}|^2d\bx+\int_{\Omega}\frac{p}{\rho}\frac{d \rho}{d c}\frac{D c}{Dt}d\bx+\int_{\partial\Omega_w}\lambda\gamma^2\rho\partial_n c \frac{D_{\Gamma} c}{Dt}dS.
	\ee
	
	Multiplying the Navier-stokes equation \eqref{nseq} by $\vel$ followed by integration by parts, the rate of change of kinetic energy is calculated as
	\be\label{kinetic}
	&&\frac{d}{dt}\int_{\Omega}\rho\frac{|\vel|^2}2d\bx\nonumber\\ &=&-\int_{\Omega}2\eta|\bD_{\eta}|^2d\bx-\int_{\Omega}\lambda|\nabla\cdot\vel|^2d\bx +\int_{\Omega}p\nabla\cdot \vel d\bx +\int_{\Omega}\lambda\gamma^2\rho(\nabla c\otimes\nabla c):\nabla\vel d\bx\nonumber\\
	&&-\int_{\partial\Omega_w}\beta_{\Gamma}|\vel^{s}|^2+\int_{\partial\Omega_w}\frac{df_w}{d c}\nabla_{\Gamma} c \cdot\vel_\tau dS\nonumber\\
	&=&-\int_{\Omega}2\eta|\bD_{\eta}|^2d\bx-\int_{\Omega}\lambda|\nabla\cdot\vel|^2d\bx +\int_{\Omega}p\nabla\cdot \vel d\bx +\int_{\Omega}\lambda\gamma^2\rho(\nabla c\otimes\nabla c):\nabla\vel d\bx\nonumber\\
	&&-\int_{\partial\Omega_w}\beta_{\Gamma}|\vel^{s}|^2-\frac{d}{dt}\int_{\partial\Omega_w}f_w dS-\int_{\partial\Omega_w}\frac{df_w}{d c}\mathcal{M}_{\Gamma}L (c) dS,
	\ee
	where we have used the definition $\frac{D_{\Gamma}c}{Dt} = \frac{\partial c}{\partial t}+\vel_{\tau}\cdot\nabla_{\Gamma} c$, $\vel_{\tau}=\vel^s$ and boundary conditions  in \eqref{model_bd}.
	
	From the derivation of $I_2$ in Appendix A, we get
	\be\label{eq32}
	&&\frac{d}{d t}\int_{\Omega}\rho\lambda_{ c}\left(G( c)+\frac{\gamma^2}{2}|\nabla c|^2\right)d\bx\nonumber\\
	&=&\int_{\Omega}\rho\lambda_{ c}\frac {DG} {D t} d\bx+ \int_{\Omega} \rho\lambda_{ c} \gamma^2\nabla c\cdot\nabla\left(\frac{D c}{D t}\right)d\bx-\int_{\Omega}\rho\lambda_{ c} \gamma^2(\nabla c \otimes\nabla c):\nabla\vel d\bx.
	\ee
	
	Combining  the equations \eqref{pheenergy}-\eqref{eq32} leads to the final energy dissipation law. $\qedsymbol$
	
	\begin{rmk}
		In the above derivation, we have neglected the external body force, for example the gravity. If the effect of gravity needs to be taken into consideration, an extra gravitational potential should be added to the total energy
		\begin{equation}
		E^{tot} =  \int_{\Omega}\frac{\rho|\vel|^2}{2} d\bx +\int_{\Omega}\lambda_{ c}\rho\left(G( c) +\frac{\gamma^2}{2}|\nabla c|^2\right)d\bx 
		+\int_{\partial\Omega_w}f_w(c)dS +\int_{\Omega}\rho g zd\bx,
		\end{equation}
		where $g$ is the gravitational constant and $z$ is the vertical position. Using the fact that \cite{guo2015thermodynamically}
		\begin{equation}
		\frac {d}{dt} \int_{\Omega}\rho g z dv = \int_{\Omega}\rho g \vel\cdot \boldsymbol{e}_zd\bx,\quad \boldsymbol{e}_z= (0,0,1)^T,
		\end{equation}
		the conservation of momentum equation \eqref{nseq} is changed to
		\begin{equation}
		\rho \frac{D\vel}{D t}  = \nabla\cdot(2\eta \bD_{\eta}) + \nabla(\lambda\nabla\cdot\vel)-\nabla p -\nabla\cdot(\lambda_{ c}\gamma^2\rho\nabla c\otimes\nabla c)-\rho g \boldsymbol{e}_z.
		\end{equation}
		
	\end{rmk}
	
	\subsection{Non-dimensionalization and Reformulation}
	In the following parts of the article,  we assume that $\lambda = -\frac{2\eta}{3}$ for simplicity. 
	Now we introduce the dimensionless variables 
	\be
	&\hat{x} = \frac{x}{L^*},\quad \hat{\vel} = \frac{\vel}{U^*},\quad  \hat\rho = \frac{\rho}{\rho*},\quad \hat{t}= \frac{t}{t^*}, \quad \hat{\tilde{\mu}} = \frac{\tilde{\mu}}{\mu^*},\quad \hat{p}= \frac{p}{P^*},\nonumber\\
	&\hat{\eta}= \frac{\eta}{\eta^*},\quad M = \frac{\mathcal{M}}{\mathcal{M}^*}, \quad M_\Gamma = \frac{\mathcal{M}_\Gamma}{ \mathcal{M}_\Gamma^*}, \quad l_s^{-1} = \frac{\beta_{\Gamma}}{\beta_{\Gamma}^*},\quad \epsilon =\frac{\gamma}{L^*}. 
	\ee
	Here $L^*$, $U^*$, $\rho^*$, $\eta^*$ and $\mathcal{M}^*$ are the characteristic scales of length, velocity, density, viscosity, and mobility coefficient, which are defined as 
	\be
	t^*=\frac{L^*}{U^*}, \quad \mu^*=\lambda_c\epsilon, \quad P^* = \rho^*\mu^*,\quad  \mathcal M^*= \frac{\rho^*U^*L^*}{\mu^*},\quad \mathcal{M}_\Gamma^* = \frac{1}{t^* \rho^* \lambda_c \gamma}, \quad \beta_{\Gamma}^{*}  = \frac{\eta^*}{L^*}.
	\ee
	For convenience, the hat symbol will be removed in the dimensionless quantities, and the dimensionless system  of (\ref{model},\ref{model_bd}) is given by
	\begin{subequations}
		\label{modelnd}
		\begin{eqnarray} 
		&&\rho \frac{Dc}{Dt} = \nabla\cdot(M\nabla\tilde{\mu}),\\
		&&\tilde{\mu} = \frac{1}{\epsilon}\frac{dG}{dc}-\frac{\epsilon}{\rho}\nabla\cdot(\rho \nabla c)+\alpha p,\\
		&&Re \rho \frac{D\vel}{Dt} = \nabla\cdot(\eta(\nabla\vel+(\nabla\vel)^T))-\nabla(\frac{2\eta}{3}\nabla\cdot\vel)-\frac{Re}{\beta}\nabla p -\frac{Re}{\beta}\epsilon\nabla\cdot(\rho\nabla c\otimes\nabla c),\\
		&&\nabla\cdot\vel =  \alpha \nabla\cdot(M\nabla\tilde{\mu}),
		\end{eqnarray}
	\end{subequations}
	with boundary conditions 
	\begin{subequations}
		\label{modelndbd}
		\begin{align} 
		\frac{D_{\Gamma}c}{Dt} = -M_{\Gamma}L(c),\\
		\partial_n\tilde{\mu}=0,\\
		\vel\cdot\bn = 0,\\
		l_s^{-1}u^s_{\tau_i} = -\bn\cdot\bsigma_{\eta}\cdot\btau_i +\frac{Re}{\beta}  L\partial_{\tau_i}c, 
		\end{align}
	\end{subequations}
	where
	\be
	L(c) = \epsilon\rho\partial_n c +\alpha_w \frac{df_w}{dc},\quad f_w( c) = -\frac{1}{2}\cos(\theta_s)\sin(\frac{(2 c-1)\pi}{2}).
	\ee
	and with dimensionless parameters
	\be
	Re = \frac{\rho^*U^*L^*}{\eta^*},\quad \beta = \frac{(U^*)^2}{\mu^*},\quad \alpha_w=\frac{\sigma}{\rho^*\lambda_c\gamma}.
	\ee
	
	
	If we define 
	\begin{subequations}
		\begin{eqnarray}
		&&\bar{p} = p+\rho(\frac{1}{\epsilon}G(c)+\frac{\epsilon}{2}|\nabla c|^2),\label{barp}\\
		&&\bar{\mu} = \tilde{\mu} -\alpha \bar{p}, \label{barmu}
		\end{eqnarray} 
	\end{subequations}
	then the system \eqref{modelnd} could be rewritten as
	\begin{subequations}
		\label{modelnd2}
		\begin{eqnarray}
		&&\rho \frac{Dc}{Dt} =  \nabla\cdot(M\nabla\bar{\mu})+ \alpha\nabla\cdot(M\nabla\bar{p}),\\
		&&	\bar{\mu} = \frac{1}{\epsilon}\frac{dG}{dc}-\frac{\epsilon}{\rho}\nabla\cdot(\rho \nabla c)-\alpha \rho(\frac{1}{\epsilon}G(c)+\frac{\epsilon}{2}|\nabla c|^2),\\
		&&	Re \rho \frac{D\vel}{Dt} = \nabla\cdot(\eta(\nabla\vel+(\nabla\vel)^T))-\nabla(\frac{2\eta}{3}\nabla\cdot\vel)-\frac{Re}{\beta}\nabla \bar{p} +\frac{Re}{\beta}\rho\bar{\mu}\nabla c\\
		&&	\nabla\cdot\vel =   \alpha\nabla\cdot(M\nabla\bar{\mu})+  \alpha^2\nabla\cdot(M\nabla\bar{p}),
		\end{eqnarray}
	\end{subequations}
	
	If we define the Sobolev spaces as follows
	\begin{eqnarray}
	\boldsymbol{W}^{1,3}(\Omega) = (W^{1,3}(\Omega))^2,\\
	\boldsymbol{W}_b^{1,3}(\Omega) = \{\vel=(u_x,u_y)^T\in\boldsymbol{W}^{1,3}|u_y=b, \mbox{on~} \partial\Omega_w\},\\
	\boldsymbol{W}_b = W^{1,3}(\Omega)\times W^{1,3}(\Omega) \times \boldsymbol{W}_b^{1,3}(\Omega)\times W^{1,3/2}(\Omega),
	\end{eqnarray}
	and 
	then above system satisfies the following energy dissipation law.
	\begin{thm}
		If $ (c, \mu, \vel, \bar{p}) \in \boldsymbol{W}_b,$ are smooth solutions of above system \eqref{modelnd2} with boundary conditions \eqref{modelndbd}, then
		the following energy law is satisfied:
		\be
		\frac{d\mathcal{E}^{tot}}{dt} &\!\!\!=\!\!\!& \frac{d}{dt}\left\{  \int_{\Omega}\frac{\rho|\vel|^2}{2} d\bx +\frac{1}{\beta}\int_{\Omega}\rho\left(\frac 1 {\epsilon}G( c) +\frac{\epsilon}{2}|\nabla c|^2\right)d\bx
		+\frac{\alpha_w }{\beta}\int_{\partial\Omega_w} f_wdS\right\}\nonumber\\
		&\!\!\!=\!\!\!&-\frac{1}{Re}\int_{\Omega}\eta\sum_{i< j}|\partial_iu_j+\partial_ju_i|^2d\bx-\frac{2}{3Re}\int_{\Omega}\eta\sum_{i<j}|\partial_iu_i-\partial_ju_j|^2d\bx  -\frac{1}{\beta}\int_{\Omega}M^n|\nabla\tilde\mu|^2d\bx\nonumber\\
		&&-\int_{\partial\Omega_w}\left(\frac{1}{\beta}M_{\Gamma}\left|L( c)\right|^2 +\frac{1}{l_sRe}|\vel^{s}|^2\right)dS.
		\ee 
	\end{thm}
	The proof  is similar as Theorem \ref{energylaw}. Here we omit the details. 
	
	\begin{rmk}
		When the walls move, i.e. $\vel_w\neq 0$, the above energy dissipation law has an extra term induced by external energy input
		\be
		\frac{d\mathcal{E}^{tot}}{dt} 
		&\!\!\!=\!\!\!&-\frac{1}{Re}\int_{\Omega}\eta\sum_{i< j}|\partial_iu_j+\partial_ju_i|^2d\bx-\frac{2}{3Re}\int_{\Omega}\eta\sum_{i<j}|\partial_iu_i-\partial_ju_j|^2d\bx  -\frac{1}{\beta}\int_{\Omega}M^n|\nabla\tilde\mu|^2d\bx\nonumber\\
		&&-\int_{\partial\Omega_w}\left(\frac{1}{\beta}M_{\Gamma}\left|L( c)\right|^2 +\frac{1}{l_sRe}|\vel^{s}|^2\right)dS-\int_{\partial\Omega_w}\frac{1}{l_sRe}\vel^{s}\cdot \vel_w dS.
		\ee  
	\end{rmk}	
	\section{Numerical Scheme and Analysis}\label{section:numerial}
	\subsection{Time-discrete primitive method}
	In this section, we present the  numerical method of system \eqref{modelnd2} with boundary conditions \eqref{modelndbd} in the primitive variable formulation. Let $\Delta t > 0$ denote the time
	step, and assume $( c^n, \bar{\mu}^n, \vel^n, \bar{p}^n)$
	are the solutions at the time $t = n\Delta t$. We then find the solutions at time $t = (n + 1)\Delta t$
	are $\left(c^{n+1}, \bar{\mu}^{n+1}, \vel^{n+1}, \bar{p}^{n+1}\right)$ that satisfy 
	
	\begin{subequations}
		\label{num1}
		\begin{eqnarray}
		&& \rho^{n}\frac{c^{n+1}-c^n}{\Delta t} + \rho^{n+1}(\vel^{n+1}\cdot\nabla)c^{n+1}	=  \nabla\cdot (M^n\nabla \bar\mu^{n+1})+  \alpha  \nabla\cdot (M^n \nabla \bar p^{n+1}),\label{numc}\\
		&& \rho^{n}\bar{\mu}^{n+1} =\frac{  \rho^{n+1/2}}{\epsilon}g(c^{n+1},c^n)-\epsilon\nabla\cdot( \rho^{n+1/2}\nabla c^{n+1/2})\nonumber\\
		&&~~~~~~~~~-\alpha \rho^n \rho^{n+1}\left(\frac{G^{n+1/2}}{\epsilon}+\frac{\epsilon }{2} (|\nabla c|^2)^{n+1/2}\right),\label{nummu}\\
		&& \rho^n \frac{\vel^{n+1}-\vel^{n}}{\Delta t}+ \rho^n(\vel^n\cdot\nabla)\vel^{n+1}+\frac 1 2 \left(\frac{ \rho^{n+1}- \rho^n}{\Delta t}+\nabla\cdot( \rho^n\vel^n))\right) +\frac{1}{\beta}\nabla \bar p^{n+1}\nonumber\\
		&&=\frac{1}{\beta} \rho^{n+1}\bar\mu^{n+1}\nabla c^{n+1}+\frac{1}{Re} \nabla\cdot(\eta^n(\nabla\vel^{n+1}+(\nabla\vel^{n+1})^T))-\frac{2}{3Re} \nabla(\eta^n\nabla\cdot\vel^{n+1}),\label{numns}\\
		&&\nabla\cdot\vel^{n+1} = \alpha   \nabla\cdot (M^n  \nabla \bar\mu^{n+1})+  \alpha^2   \nabla\cdot (M^n \nabla \bar p^{n+1}),\label{numdiv} 
		\end{eqnarray}
	\end{subequations}
	with boundary conditions 
	\begin{subequations}
		\label{num1bd}
		\begin{eqnarray}	&&\frac{c^{n+1}-c^n}{\Delta t} +\vel_{\btau}^{n+1}\cdot\nabla_{\Gamma} c^{n+1/2}= -M_{\Gamma}L^{n+1/2}(c),\label{numcbd}\\
		&&\partial_n\tilde{\mu}^{n+1}=0,\label{nummubd}\\
		&& L^{n+1/2}(c) = \epsilon \rho^{n+1/2}\partial_n c^{n+1/2} +\alpha_w \frac{f_w(c^{n+1})-f_w(c^n)}{c^{n+1}-c^n},\label{numl}\\
		&&\vel^{n+1}\cdot\bn = 0,\\
		&&l_s^{-1}u_{\tau_i}^{s,n+1} = -\bn\cdot(\eta^n(\nabla\vel^{n+1}+(\nabla\vel^{n+1})^T))\cdot\btau_i +\frac{Re}{\beta}  L^{n+1/2}\partial_{\tau_i}c^{n+1/2}, \label{numslip}
		\end{eqnarray}
	\end{subequations}
	where $\tilde{\mu}$ is defined in \eqref{barmu} and we have used the notations 
	\be
	& (\cdot)^{n+1/2}=\frac 1 2 [(\cdot)^{n+1}+(\cdot)^n], \quad \rho^{n+1} =  \rho(c^{n+1}),\\
	& g(c^{n+1},c^n) = \frac 1 4 (c^{n+1}(c^{n+1}-1)+c^n(c^n-1))(c^{n+1}+c^n-1).
	\ee
	
	For above discretization, it satisfies the following properties. 
	\begin{lma}(\cite{guo_mass_2017})
		If $c^{n+1}$ is the solution of above system \eqref{num1}-\eqref{num1bd}, then we have 
		\begin{subequations}
			\begin{align}
			G(c^{n+1})-G(c^n) = g(c^{n+1},c^{n})(c^{n+1}-c^{n}),\label{G_c}\\
			\rho(c^{n+1})-\rho(c^n) = -\alpha  \rho^{n+1} \rho^{n}(c^{n+1}-c^{n}).\label{rho_c}
			\end{align}
		\end{subequations}
		And the system \eqref{num1}-\eqref{num1bd} yields mass conservation for each component of binary fluid
		\begin{subequations}
			\begin{align}
			\int_{\Omega} \rho^{n+1}d\bx = \int_{\Omega} \rho^{n}d\bx,\\
			\int_{\Omega} \rho^{n+1}c^{n+1}d\bx = \int_{\Omega} \rho^{n}c^n d\bx.
			\end{align}
		\end{subequations}	
		
	\end{lma}

	\begin{thm}\label{numstable}
		If $ \left(c^{n+1}, \mu^{n+1}, \vel^{n+1}, \bar{p}^{n+1}\right)\in \boldsymbol{W}_b $ are  solutions of above system \eqref{num1} with boundary conditions \eqref{num1bd}, then
		the following energy law is satisfied:
		\be
		&&\mathcal{E}^{n+1,tot} -\mathcal{E}^{n,tot}\nonumber\\
		&=&-\frac{\Delta t}{Re}\int_{\Omega}\eta^n\sum_{i< j}|\partial_iu^{n+1}_j+\partial_ju^{n+1}_i|^2d\bx-\frac{2\Delta t}{3Re}\int_{\Omega}\eta^{n}\sum_{i<j}|\partial_iu^{n+1}_i-\partial_ju^{n+1}_j|^2d\bx \nonumber\\
		&& -\frac{\Delta t  }{\beta  }\int_{\Omega}M^n|\nabla\tilde\mu^{n+1}|^2d\bx-\Delta t\int_{\partial\Omega_w}\left(\frac{1}{\beta}M_{\Gamma}\left|L^{n+1/2}( c)\right|^2 +\frac{1}{l_sRe}|\vel^{s,n+1}|^2\right)dS\\
		&&	- \Delta t\int_{\partial\Omega_w}\frac{ 1}{l_s Re}\vel^{s,n+1}\cdot \vel_w dS,
		\ee 
		where 
		$$
		\mathcal{E}^{n+1,tot} = \int_{\Omega}\frac{ \rho^{n+1}}{2} |\vel^{n+1}|^2d\bx +\frac{1}{\beta}\int_{\Omega} \rho^{n+1}\left(\frac 1 {\epsilon}G(c^{n+1}) +\frac{\epsilon}{2}|\nabla c^{n+1}|^2\right)d\bx
		+\frac{\alpha_w}{\beta}\int_{\partial\Omega_w} f_w(c^{n+1})dS
		$$ 
		is the discretized total energy. 
	\end{thm}
	\noindent \textbf{Proof:} Taking inner product of the first equation \eqref{numc} with $\frac{\Delta t}{\beta}\bar\mu^{n+1}$ results in the following equation
	\be\label{numenergyc}
	&&\frac{1}{\beta}\int_{\Omega} \rho^{n}(c^{n+1}-c^n)\bar\mu^{n+1}d\bx +\frac{\Delta  t}{\beta}\int_{\Omega} \rho^{n+1}\vel^{n+1}\cdot\nabla c^{n+1} \bar\mu^{n+1}d\bx\nonumber\\
	&=&-\frac{\Delta t  }{\beta  }\int_{\Omega}M^n\nabla\bar\mu^{n+1}\cdot \nabla\tilde\mu^{n+1} d\bx,
	\ee
	where we used the boundary condition \eqref{nummubd} and the definition of $\bar\mu$ in \eqref{barmu}.
	
	Multiplying the second equation \eqref{nummu} with $\frac{c^{n+1}-c^n}{\beta}$ yields \cite{guo_mass_2017}
	\be\label{numenergymu}
	&&\frac{1}{\beta}\int_{\Omega} \rho^{n}(c^{n+1}-c^n)\bar\mu^{n+1}d\bx \nonumber\\
	&=& \frac{1}{\beta}\int_{\Omega} \rho^{n+1}\left(\frac{1}{\epsilon}G(c^{n+1})+\frac{\epsilon} 2  |\nabla c^{n+1}|^2\right)d\bx - \frac{1}{\beta}\int_{\Omega} \rho^{n}\left(\frac{1}{\epsilon}G(c^{n})+\frac{\epsilon}2 |\nabla c^{n}|^2\right)d\bx\nonumber\\ &&-\frac{1}{\beta}\int_{\partial\Omega_w}\epsilon \rho^{n+1/2}\partial_nc^{n+1/2}(c^{n+1}-c^n)dS,
	\ee
	where we use the results in above Lemma and the boundary condition \eqref{numcbd}.
	
	Multiplying the Navier-Stokes equation \eqref{numns} with $\Delta t\vel^{n+1}$, we have
	\be\label{numenergyns}
	&&\frac{1}{2}\int_{\Omega}( \rho^{n+1}|\vel^{n+1}|^2)d\bx -\frac{1}{2}\int_{\Omega}( \rho^{n}|\vel^{n}|^2)d\bx\nonumber\\
	&=&-\frac{\Delta t}{Re}\int_{\Omega}\eta^n\sum_{i< j}|\partial_iu^{n+1}_j+\partial_ju^{n+1}_i|^2d\bx-\frac{2\Delta t}{3Re}\int_{\Omega}\eta^{n}\sum_{i<j}|\partial_iu^{n+1}_i-\partial_ju^{n+1}_j|^2d\bx \nonumber\\
	&&+\frac{\Delta t}{\beta}\int_{\Omega}\bar p^{n+1}\nabla\cdot \vel^{n+1}d\bx +\frac{\Delta t}{\beta}\int_{\Omega} \rho^{n+1}\bar\mu^{n+1}\nabla c^{n+1}\cdot\vel^{n+1}d\bx\nonumber\\
	&& - \frac{\Delta t}{Re}\int_{\partial\Omega_w}l_s^{-1}|\vel^{s,n+1}|^2dS+ \frac{\Delta t}{\beta}\int_{\partial\Omega_w}  L^{n+1/2}\nabla_{\Gamma}c^{n+1/2}\cdot \vel_{\btau}^{n+1}dS,\nonumber\\
	&&-\frac{\Delta t}{Re}\int_{\partial\Omega_w} l_s^{-1}\vel^{s,n+1}\cdot \vel_w dS,
	\ee 
	where we used the slip boundary condition \eqref{numslip} and the tensor calculation in Appendix \ref{tensor_calculation}.
	
	For the last term in above equation, combining the definition of $L^{n+1/2}$ in \eqref{numl} and equation \eqref{numcbd} yields
	\be\label{numenergsp}
	&&\frac{\Delta t}{\beta}\int_{\partial\Omega_w}  L^{n+1/2}\nabla_{\Gamma}c^{n+1/2}\cdot \vel_{\btau}^{n+1}dS\nonumber\\
	&=&\frac{\Delta t}{\beta}\left(\int_{\partial\Omega_w} \alpha_w\frac{f_w(c^{n+1})-f_w(c^n)}{c^{n+1}-c^n} \nabla_{\Gamma}c^{n+1/2}\cdot \vel_{\btau}^{n+1}dS\right.\nonumber\\
	&&\left.+\int_{\partial_\Omega}\epsilon \rho^{n+1/2}\partial_nc^{n+1/2}\nabla_{\Gamma}c^{n+1/2}\cdot \vel_{\btau}^{n+1}dS\right)\nonumber\\
	&=&-\frac{1}{\beta}\int_{\partial\Omega_w}\alpha_w (f_w(c^{n+1})-f_w(c^n))dS-\frac{\Delta t}{\beta}\int_{\partial\Omega_w}\alpha_w \frac{f_w(c^{n+1})-f_w(c^n)}{c^{n+1}-c^n}M_{\Gamma}L^{n+1/2}dS\nonumber\\
	&&+\frac{\Delta t}{\beta}\int_{\partial_\Omega}\epsilon \rho^{n+1/2}\partial_nc^{n+1/2}\nabla_{\Gamma}c^{n+1/2}\cdot \vel_{\btau}^{n+1}dS.
	\ee 
	
	Then equation \eqref{numenergyns} could be rewritten as 
	\be\label{numkinetic}
	&&\frac{1}{2}\int_{\Omega}( \rho^{n+1}|\vel^{n+1}|^2)d\bx -\frac{1}{2}\int_{\Omega}( \rho^{n}|\vel^{n}|^2)d\bx
	+\frac{1}{\beta}\int_{\partial\Omega_w}\alpha_w (f_w(c^{n+1})-f_w(c^n))dS\nonumber\\
	&=&-\frac{\Delta t}{Re}\int_{\Omega}\eta^n\sum_{i< j}|\partial_iu^{n+1}_j+\partial_ju^{n+1}_i|^2d\bx-\frac{2\Delta t}{3Re}\int_{\Omega}\eta^{n}\sum_{i<j}|\partial_iu^{n+1}_i-\partial_ju^{n+1}_j|^2d\bx \nonumber\\
	&&+\frac{\Delta t}{\beta}\int_{\Omega}\bar p^{n+1}\nabla\cdot \vel^{n+1}d\bx +\frac{\Delta t}{\beta}\int_{\Omega} \rho^{n+1}\bar\mu^{n+1}\nabla c^{n+1}\cdot\vel^{n+1}d\bx\nonumber\\
	&&- \frac{\Delta t}{Re}\int_{\partial\Omega_w}l_s^{-1}|\vel^{s,n+1}|^2dS -\frac{\Delta t}{\beta}\int_{\partial\Omega_w}\alpha_w \frac{f_w(c^{n+1})-f_w(c^n)}{c^{n+1}-c^n}M_{\Gamma}L^{n+1/2}dS\nonumber\\
	&&+\frac{\Delta t}{\beta}\int_{\partial_\Omega}\epsilon \rho^{n+1/2}\partial_nc^{n+1/2}\nabla_{\Gamma}c^{n+1/2}\cdot \vel_{\btau}^{n+1}dS-\frac{\Delta t}{Re}\int_{\partial\Omega_w} l_s^{-1}\vel^{s,n+1}\cdot \vel_w dS.
	\ee
	
	Multiplying the last equation \eqref{numdiv} with $\frac{\Delta t}{\beta}\bar p^{n+1}$ yields
	\be\label{numenergydiv}
	0= -\int_{\Omega}\frac{\Delta t}{\beta}\bar p^{n+1}\nabla \cdot \vel^{n+1}d\bx -\frac{\Delta t  }{ \beta}\int_{\Omega}M^n\nabla\tilde\mu^{n+1}\cdot\nabla\bar p^{n+1}d\bx.
	\ee
	
	Summing up equations \eqref{numenergyc}, \eqref{numenergymu}, \eqref{numkinetic} and \eqref{numenergydiv} results 
	\be
	&&\mathcal{E}^{n+1,tot} -\mathcal{E}^{n,tot}\nonumber\\
	&\!\!\!=\!\!\!&-\frac{\Delta t}{Re}\int_{\Omega}\eta^n\sum_{i< j}|\partial_iu^{n+1}_j+\partial_ju^{n+1}_i|^2d\bx-\frac{2\Delta t}{3Re}\int_{\Omega}\eta^{n}\sum_{i<j}|\partial_iu^{n+1}_i-\partial_ju^{n+1}_j|^2d\bx \nonumber\\
	&& -\frac{\Delta t  }{\beta  }\int_{\Omega}M^n|\nabla\tilde\mu^{n+1}|^2d\bx-\int_{\partial\Omega_w}\frac{\Delta t}{l_sRe}|\vel^{s,n+1}|^2dS -\int_{\partial\Omega_w} \frac{\Delta t}{l_sRe}\vel^{s,n+1}\cdot \vel_w dS\nonumber\\
	&& -\frac{\Delta t}{\beta}\int_{\partial\Omega_w}\alpha_w \frac{f_w(c^{n+1})-f_w(c^n)}{c^{n+1}-c^n}M_{\Gamma}L^{n+1/2}dS\nonumber\\
	&&+\frac{1}{\beta}\int_{\partial\Omega_w}\epsilon \rho^{n+1/2}\partial_nc^{n+1/2}(c^{n+1}-c^n)dS +\frac{\Delta t}{\beta}\int_{\partial_\Omega}\epsilon \rho^{n+1/2}\partial_nc^{n+1/2}\nabla_{\Gamma}c^{n+1/2}\cdot \vel_{\tau}^{n+1}dS\nonumber\\
	&\!\!\!=\!\!\!&-\frac{\Delta t}{Re}\int_{\Omega}\eta^n\sum_{i< j}|\partial_iu^{n+1}_j+\partial_ju^{n+1}_i|^2d\bx-\frac{2\Delta t}{3Re}\int_{\Omega}\eta^{n}\sum_{i<j}|\partial_iu^{n+1}_i-\partial_ju^{n+1}_j|^2d\bx \nonumber\\
	&& -\frac{\Delta t  }{\beta  }\int_{\Omega}M^n|\nabla\tilde\mu^{n+1}|^2d\bx-\int_{\partial\Omega_w}\frac{\Delta t}{l_sRe}|\vel^{s,n+1}|^2dS-\int_{\partial\Omega_w} \frac{\Delta t}{l_sRe}\vel^{s,n+1}\cdot \vel_w dS\nonumber\\
	&& -\frac{\Delta t}{\beta}\int_{\partial\Omega_w}\alpha_w \frac{f_w(c^{n+1})-f_w(c^n)}{c^{n+1}-c^n}M_{\Gamma}L^{n+1/2}dS\nonumber-\frac{\Delta t}{\beta}\int_{\partial\Omega_w}\epsilon \rho^{n+1/2}\partial_nc^{n+1/2}M_{\Gamma}L^{n+1/2}dS\nonumber\\
	&\!\!\!=\!\!\!&-\frac{\Delta t}{Re}\int_{\Omega}\eta^n\sum_{i< j}|\partial_iu^{n+1}_j+\partial_ju^{n+1}_i|^2d\bx-\frac{2\Delta t}{3Re}\int_{\Omega}\eta^{n}\sum_{i<j}|\partial_iu^{n+1}_i-\partial_ju^{n+1}_j|^2d\bx \nonumber\\
	&& -\frac{\Delta t  }{\beta  }\int_{\Omega}M^n|\nabla\tilde\mu^{n+1}|^2d\bx-\int_{\partial\Omega_w}\frac{\Delta t}{l_sRe}|\vel^{s,n+1}|^2dS -\frac{\Delta t}{\beta}\int_{\partial\Omega_w}M_{\Gamma}|L^{n+1/2}|^2dS\nonumber\\
	&&-\int_{\partial\Omega_w} \frac{\Delta t}{l_sRe}\vel^{s,n+1}\cdot \vel_w dS,
	\ee
	where we used the definition of $\bar\mu$ \eqref{barmu} and slip boundary condition \eqref{numslip}.  $\qedsymbol$

	\subsection{Fully-discrete $C^0$ finite element scheme}
	The fully-discrete  $C^0$ finite element scheme for this time-discrete primitive scheme  \eqref{num1}-\eqref{num1bd} is presented in the section. For simplicity, we  only consider a two-dimensional case here. It is straightforward to extend the results to three-dimensional case.
	The domain $\Omega$ is a bounded domain with Lipschitz-continuous boundary $\partial\Omega$. Specifically, we denote $\partial\Omega_w$ as the solid wall where the slip boundary condition is used.   
	Let $\boldsymbol{W}^h_b = H^h\times H^h \times \boldsymbol{U}^h_b\times P^h$
	be the finite dimensional space of $\boldsymbol{W}_b$ based on a given finite element discretization of $\Omega$. 
	If we assume that $ \rho^{n}\in L^{\infty}(\Omega)$ and positive \cite{guo2014numerical}, then the weak form of semi-discrete system  \eqref{num1} with boundary conditions \eqref{num1bd} is the following: finding 
	$(c_h^{n+1}, \bar{\mu}_h^{n+1}, \vel_h^{n+1}, \bar{p}_h^{n+1})\in \boldsymbol{W}_b^h$, such that 
	\begin{subequations}
		\label{num_fem}
		\begin{eqnarray}
		&&\int_{\Omega}\left( \rho_h^{n}\frac{c_h^{n+1}-c_h^n}{\Delta t} + \rho_h^{n+1}(\vel_h^{n+1}\cdot\nabla)c_h^{n+1}\right)\psi_h d\bx\nonumber\\
		&&= - \int_{\Omega} M^n  \left(\nabla\bar\mu_h^{n+1}+\alpha \nabla \bar p_h^{n+1}\right)\cdot\nabla\psi_h d\bx, \label{numc_fem}\\
		&&	\int_{\Omega} \rho_h^{n}\bar{\mu}_h^{n+1}\chi_h d\bx =\int_{\Omega}\frac{  \rho_h^{n+1/2}}{\epsilon}g(c_h^{n+1},c_h^n)\chi_h d\bx+\int_{\Omega}\epsilon \rho_h^{n+1/2}\nabla c_h^{n+1/2}\cdot\nabla\chi_h d\bx \nonumber\\
		&&-\int_{\Omega}\alpha \rho_h^n \rho_h^{n+1}\left(\frac{G_h^{n+1/2}}{\epsilon}+\frac{\epsilon }{2} (|\nabla c_h|^2)^{n+1/2}\right)\chi_h d\bx -\int_{\partial\Omega_w}\epsilon \rho_h^{n+1/2}\partial_n c_h^{n+1/2}\chi_h dS,\label{nummu_fem}\\
		&&\int_{\Omega}\left\{ \rho_h^n \frac{\vel_h^{n+1}-\vel_h^{n}}{\Delta t}+ \rho_h^n(\vel_h^n\cdot\nabla)\vel_h^{n+1}+\frac 1 2 \left(\frac{ \rho_h^{n+1}- \rho_h^n}{\Delta t}+\nabla\cdot( \rho_h^n\vel_h^n)\right)\right\}\cdot\bv_h d\bx \nonumber\\
		&&=-\frac{1}{\beta}\int_{\Omega}\nabla \bar p_h^{n+1}\cdot \bv_h d\bx\int_{\Omega}\frac{1}{\beta} \rho_h^{n+1}\bar\mu_h^{n+1}\nabla c^{n+1}\cdot \bv_h d\bx\nonumber\\
		&&-\frac{1}{Re}\int_{\Omega} (\eta^n(\nabla\vel_h^{n+1}-(\nabla\vel_h^{n+1})^T)):\nabla\bv_h d\bx+\frac{2}{3Re}\int_{\Omega}\eta^n\nabla\cdot\vel_h^{n+1}\nabla\cdot\bv_hd\bx\nonumber\\ &&-\int_{\partial\Omega_w}\frac{1}{Re l_s}\vel_{h}^{s,n+1}\cdot\bv_h dS+\int_{\partial_\Omega}\frac{1}{\beta}L_h^{n+1/2}\nabla_{\Gamma}c_h^{n+1/2}\cdot\bv_h dS,\label{numns_fem}\\
		&&-\int_{\Omega}\nabla q_h\cdot\vel_h^{n+1} = - \alpha  \int_{\Omega}M^n\left( \nabla \bar\mu_h^{n+1}+ \alpha\nabla\bar p_h^{n+1}\right)\cdot\nabla q_h d\bx,\label{numdiv_fem}
		\end{eqnarray}
	\end{subequations} 
	for any $(\psi_h,\chi_h, \bv_h,q_h)\in \boldsymbol{W}_0^h$. 
	
	\begin{lma}\label{lma:mass_con}
		The fully discretized system \eqref{num_fem} satisfies mass conservation for each component of binary fluid
		\begin{subequations}
			\begin{align}
			\int_{\Omega} \rho_h^{n+1}d\bx = \int_{\Omega} \rho_h^{n}d\bx,\\
			\int_{\Omega} \rho_h^{n+1}c_h^{n+1}d\bx = \int_{\Omega} \rho_h^{n}c_h^n d\bx.
			\end{align}
		\end{subequations}	
	\end{lma}
	\noindent\textbf{Proof:} Setting $\psi_h =q_h= \rho_h^{n+1}$ in Eqs.\eqref{numc_fem} and \eqref{numdiv_fem},  we have  
	\begin{eqnarray}
	&&\int_{\Omega}\left( \rho^{n}_h \rho_h^{n+1}\frac{c_h^{n+1}-c_h^n}{\Delta t} +( \rho_h^{n+1})^2(\vel_h^{n+1}\cdot\nabla)c_h^{n+1}\right) d\bx\nonumber\\
	&&=- \int_{\Omega} M^n  \left(\nabla\bar\mu_h^{n+1}+\alpha \nabla \bar p_h^{n+1}\right)\cdot\nabla \rho^{n+1}_h d\bx,\nonumber\\
	&&-\frac{1}{\alpha}\int_{\Omega}\nabla \rho^{n+1}_h\cdot\vel_h^{n+1} d\bx=  \int_{\Omega}M^n\left( \nabla \bar\mu_h^{n+1}+ \alpha\nabla\bar p_h^{n+1}\right)\cdot\nabla \rho^{n+1}_h d\bx\nonumber.
	\end{eqnarray}
	Adding the above two equations and using Eqs. \eqref{quasi_incom2} and \eqref{rho_c} yields 
	\begin{equation}
	\int_{\Omega}\left(\frac{ \rho^{n+1}_h- \rho^{n}_h}{\Delta t} + \nabla\cdot( \rho^{n+1}\vel_h^{n+1})\right) d\bx = 0.
	\end{equation}
	Using the boundary condition $\vel^{n+1}_h \cdot\bn = 0$, we have the conservation of total mass
	\begin{eqnarray}
	\int_{\Omega}( \rho^{n+1}_h- \rho^{n}_h)d\bx = 0.
	\end{eqnarray}

	Choosing $\psi_h =q_h=\hat \rho_h^{n+1}c^{n+1}_h$ in Eqs.\eqref{numc_fem} and \eqref{numdiv_fem},  similarly we have 
	\begin{equation}
	\int_{\Omega}\left(\frac{ \rho^{n+1}_h- \rho^{n}_h}{\Delta t}c^{n+1}_h + \nabla\cdot( \rho^{n+1}\vel_h^{n+1})\right) c^{n+1}_hd\bx = 0.
	\end{equation}
	Choosing $\psi_h =1$ in Eq.\eqref{numc_fem}  yields
	\begin{equation}
	\int_{\Omega}\left( \rho_h^{n}\frac{c_h^{n+1}-c_h^n}{\Delta t} + \rho_h^{n+1}(\vel_h^{n+1}\cdot\nabla)c_h^{n+1}\right) d\bx= 0.
	\end{equation}
	Adding the above two equations and using the velocity boundary condition $\vel^{n+1}_h \cdot\bn = 0$, we have
	\begin{equation}
	\int_{\Omega}( \rho^{n+1}_hc^{n+1}_h -  \rho^{n}_hc^{n}_h) d\bx= 0. \nonumber~~~~~~\qedsymbol
	\end{equation}  
	
	\begin{thm}\label{numstable_fem}
		If $ (c_h^{n+1}, \bar{\mu}_h^{n+1},\vel_h^{n+1}, \bar{p}_h^{n+1})$ are  solutions of the above system \eqref{num_fem}, then
		the following energy law is satisfied:
		\be
		&&\mathcal{E}_h^{n+1,tot} -\mathcal{E}_h^{n,tot}\nonumber\\
		&\!\!\!=\!\!\!&-\frac{\Delta t}{Re}\int_{\Omega}\eta_h^n\sum_{i< j}|\partial_iu^{n+1}_{h,j}+\partial_ju^{n+1}_{h,i}|^2d\bx-\frac{2\Delta t}{3Re}\int_{\Omega}\eta^{n}\sum_{i<j}|\partial_iu^{n+1}_{h,i}-\partial_ju^{n+1}_{h,j}|^2d\bx \nonumber\\
		&& -\frac{\Delta t  }{\beta  }\int_{\Omega}M^n|\nabla\tilde\mu_h^{n+1}|^2d\bx-\Delta t\int_{\partial\Omega_w}\left(\frac{1}{\beta}M_{\Gamma}\left|L_h^{n+1/2}( c_h)\right|^2 +\frac{1}{l_sRe}|\vel^{s,n+1}_{h}|^2\right)dS\nonumber\\
		&&-\Delta t\int_{\partial\Omega_w} \frac{1}{l_s Re } \vel_h^{s,n+1}\cdot \vel_w dS
		\ee 
		where 
		$$\mathcal{E}_h^{n+1} = \int_{\Omega}\frac{ \rho_h^{n+1}|\vel_h^{n+1}|^2}{2} d\bx +\frac{1}{\beta}\int_{\Omega} \rho_h^{n+1}\left(\frac 1 {\epsilon}G(c_h^{n+1}) +\frac{\epsilon}{2}|\nabla c_h^{n+1}|^2\right)d\bx
		+\frac{\alpha_w}{\beta}\int_{\partial\Omega_w} f_w(c_h^{n+1})dS$$ 
		is the discretized total energy. 
	\end{thm}
	It can be proved  by choosing 
	\begin{eqnarray}
	\psi_h = \frac{\Delta t }{\beta}\bar\mu_h^{n+1},~~
	\chi_h = \frac{c_h^{n+1}-c_h^{n}}{\beta},~~
	\bv_h = \Delta t\vel_h^{n+1},~~
	q_h = \frac{\Delta t}{\beta} \bar p^{n+1}_h
	\end{eqnarray}
	in Eq. \eqref{num_fem} and following the proof of Theorem \ref{numstable}.

	\section{Simulation Results}\label{section:simulation}
	In this section, we present some numerical simulations using the aforementioned algorithm.  Three cases are considered: Couette flow,  moving droplets in shear flow and rising bubbles to illustrate the convergence rate, the effect of  contact angle and the quasi-incompressibility of two-phase flow with large density ratio, respectively.  All of the numerical simulations in this part are based on the proposed finite element scheme and implemented with the  FreeFem++ \cite{MR3043640}. 
	

	\subsection{Convergence Study: Couette Flow }
	
\begin{figure}[!ht]
		\centering
		\includegraphics[width=7.in]{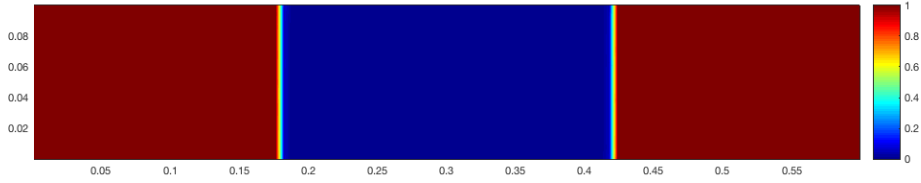} 
		\caption{Initial condition of phase 1 concentration   $c$ for the two phase Couette flow. Positions a and b are two contact points on the bottom walls. }
		\label{fig:inital}
	\end{figure}

	We start with convergence test using Couette flow with  different density and viscosity \cite{gao_efficient_2014,yu_numerical_2017} as in  Fig. \ref{fig:inital}.
	The domain size is  $ [0, 0.6]\times[0, 0.1]$ 
	The top and bottom walls move oppositely with $\mathbf{u}_w = (1,0)^T$. We do the convergence study for two-phase fluids with  both low and high  density ratios. 
	
	For the case of low density ratio, the parameters are listed as follows:
	$$Re = 200, \beta = 1.76\times 10^{-2}, M = 1.5\times 10^{-8},    \epsilon = 0.01,\alpha_w = 8.33\times10^{-4},\rho_1 = 0.8, \rho_2=1,$$
	$$ \eta_1=\eta_2 =1, \theta_s = 120^{\circ}, M_{\Gamma} = 5\times 10^5,$$
	$$ l_{s1}=l_{s2} = 0.02,  \Delta t = 8\times 10^{-4}.$$ 
	
	We first present the convergence study for P1 element with $h = 1/160, 1/226, 1/320, 1/640$ and P2 element with $h= 1/80, 1/113, 1/160, 1/320$.
	The  results  with  $ h=1/640$  and $1/320$ are used as the reference solutions for P1 and P2 elements, respectively. 

	\begin{table}[ht]
		\centering
		\begin{tabular}{|c|c|c|c|c|c|c|}
			\hline
			{Space step
				$h$} & \multicolumn{6}{|c|}{P1 Element} \\
			\cline{2-7}
			& Err($u_x$) & Rate &  Err($u_y$) & Rate& Err(c) & Rate \\ 
			\hline
			1/160& 4.9e-2 & & 2.6e-2 & & 8.5e-2 & \\ 
			\hline
			1/226& 1.9e-2 & 2.57 & 1.8e-2 & 1.46 & 5.3e-2 & 1.59 \\ 
			\hline
			1/320& 7.2e-3 & 2.64  & 6.4e-3 & 2.78 & 2.2e-2 & 2.38 \\ 
			\hline
			{Space step $h$} & \multicolumn{6}{|c|}{P2  Element} \\
			\cline{2-7}
			& Err($u_x$) & Rate &  Err($u_y$) & Rate& Err(c) & Rate \\ 
			\hline
			1/80& 2.9e-2 & & 2.5e-2 & & 7.4e-2 & \\ 
			\hline
			1/113& 9.9e-3 & 2.93 & 7.5e-3 & 3.32 & 2.8e-2 & 2.63 \\ 
			\hline
			1/160& 3.0e-3 & 3.34 & 2.5e-3 & 2.97 & 1.0e-2 & 2.80 \\ 
			\hline
		\end{tabular}
		\caption{$L^2$ norm of the error and convergence rate for  velocity $\vel=(u_x,u_y)$, phase function $c$,  at time $t=0.2 $ with   density ratio $\rho_1:\rho_2 = 0.8:1$ viscosity ratio $\eta_1:\eta_2 = 1:1$. }
		\label{tab:convergence_low}
	\end{table}
	
	
	For the case of high density ratio, the parameters are listed as follows:
	$$Re = 20, \beta = 1.76\times 10^{-2}, M = 1.5\times 10^{-8},  \rho_1 = 0.1, \rho_2=10, $$
	$$\eta_1=0.1, \eta_2 =10, \epsilon = 0.01,\alpha_w = 8.33\times 10^{-4},\theta_s = 120^{\circ}, M_{\Gamma} = 5\times 10^5,$$ 
	$$ l_{s1} = 0.01,l_{s2} = 0.0027,   \Delta t = 8\times 10^{-4}.$$ 
	The convergence rate for both P1   element and P2  element are shown in Table \ref{tab:convergence_high}. It illustrates the 2nd-order for P1 element and 3rd-order for P2 element  convergence rate in the sense of $L^2$ norm.
	
	\begin{table}[ht]
		\centering
		\begin{tabular}{|c|c|c|c|c|c|c|}
			\hline
			{Space step
				$h$} & \multicolumn{6}{|c|}{P1  Element} \\
			\cline{2-7}
			& Err($u_x$) & Rate &  Err($u_y$) & Rate& Err(c) & Rate \\ 
			\hline
			1/160& 8.8e-3& & 4.9e-3 & &2.2e-2 & \\ 
			\hline
			1/226&6.9e-3 &1.27 & 3.8e-3&1.27 & 2.0e-2& 1.12 \\ 
			\hline
			1/320& 2.2e-3& 3.17& 1.5e-3& 2.59 & 7.8e-3& 2.52 \\ 
			\hline
			{Space step $h$} & \multicolumn{6}{|c|}{P2  Element} \\
			\cline{2-7}
			& Err($u_x$) & Rate &  Err($u_y$) & Rate& Err(c) & Rate \\ 
			\hline
			1/80&5.8e-3 & &3.6e-3& &2.5e-2 & \\ 
			\hline
			1/113& 4.0e-3&1.46 &1.9e-3 &1.86 &1.1e-2 & 2.20\\ 
			\hline
			1/160& 1.1e-3& 3.83& 6.3e-4& 3.08& 3.5e-3&3.27 \\ 
			\hline
		\end{tabular}
		\caption{$L^2$ norm of the error and convergence rate for  velocity $\vel=(u_x,u_y)$, phase function $c$,  at time $T=0.2$ with   density ratio $\rho_1:\rho_2 = 0.1:10$ viscosity ratio $\eta_1:\eta_2 = 0.1:10$. }
		\label{tab:convergence_high}
	\end{table}

\begin{figure}[ht]
		\centering
		\includegraphics[width=7in]{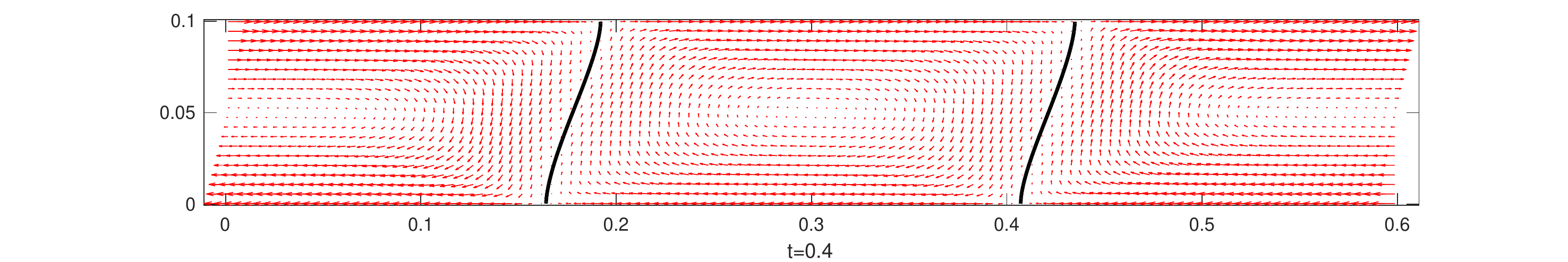}
		\caption{The interface and velocity profile at  $T=0.4$ large  density ratio $\rho_1=0.1,~\rho_2=10$. }
		\label{fig:high_interface}
	\end{figure}

	The profile of interface and velocity fields around steady state are shown in Fig. \ref{fig:high_interface}. 
	
\begin{figure}[ht]
		\centering
		\includegraphics[width=3.in]{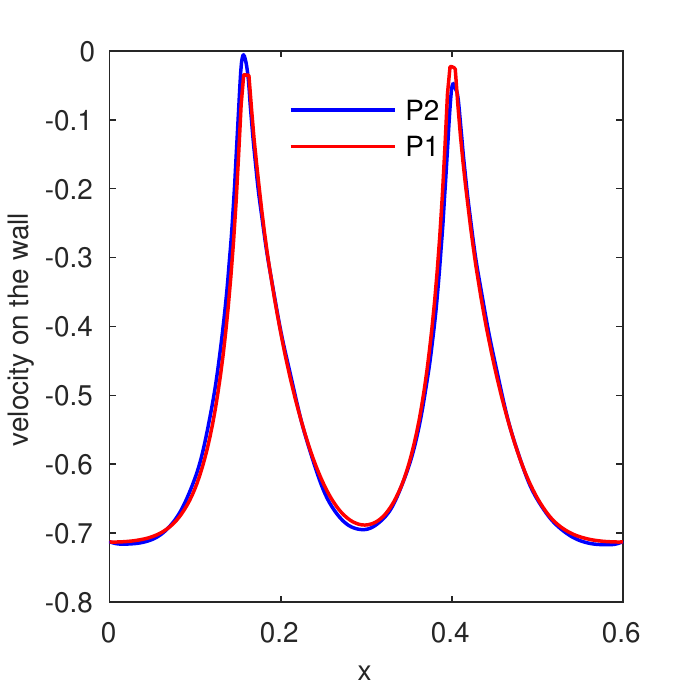}
		\includegraphics[width=3.in]{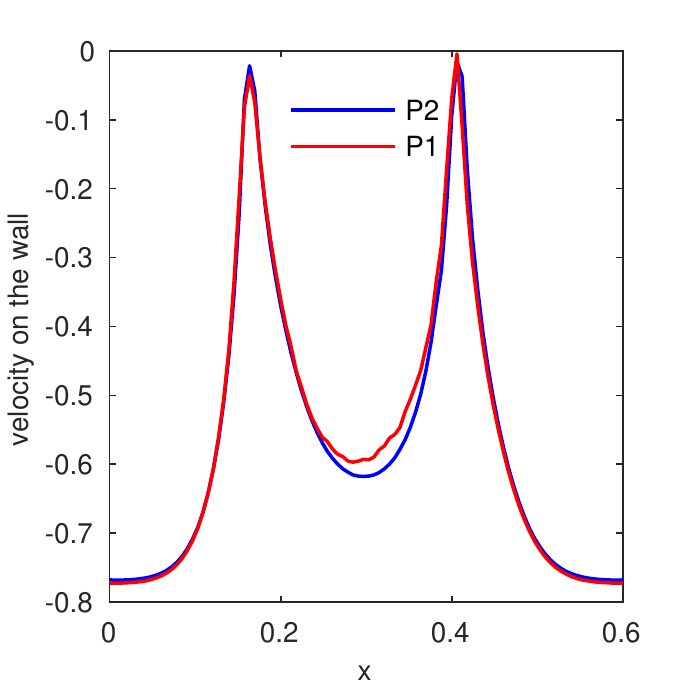}
		\caption{Velocity on wall around equilibrium state.Left: low density ratio $\rho_1 = 0.8$, $\rho_2 =1$; Right: large  density ratio $\rho_1=0.1,~\rho_2=10$. }
		\label{fig:highmcl}
	\end{figure}

	The fluid velocities on the wall are shown in Fig.\ref{fig:highmcl}. It shows that for both low and high density ratio case, P1 element could yield consistent contact velocity with P2 element.

\begin{figure}[!ht]
		\centering
		\includegraphics[width=3.5in]{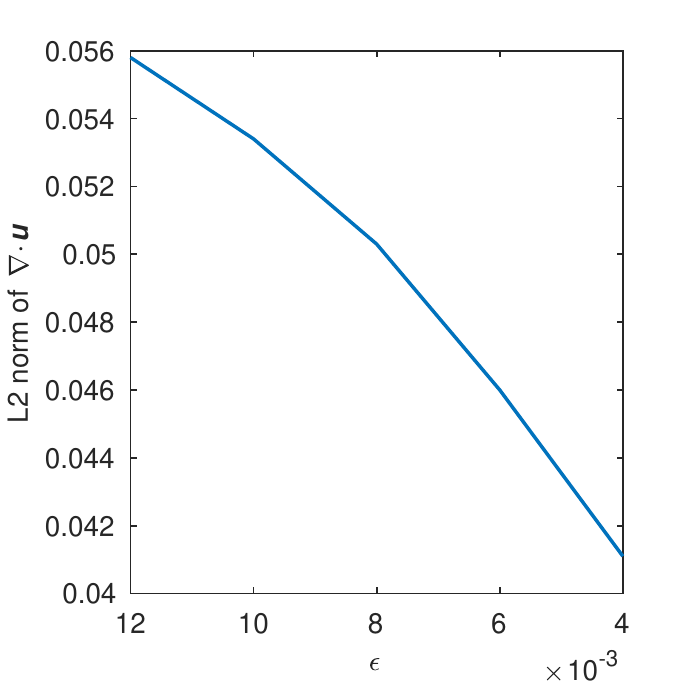}
		\caption{$L^2$ norm of $\nabla\cdot\vel$ with different $\epsilon$.}
		\label{fig:diffeps}
	\end{figure}

	In Fig.\ref{fig:diffeps}, we check the $L^2$ norm of $\nabla\cdot\vel$ with different $\epsilon$. The results confirm that as $\epsilon$ decreases, the solution converges to the sharp interface incompressible fluids.

\begin{figure}[!ht]
		\centering
		\includegraphics[width=3.in]{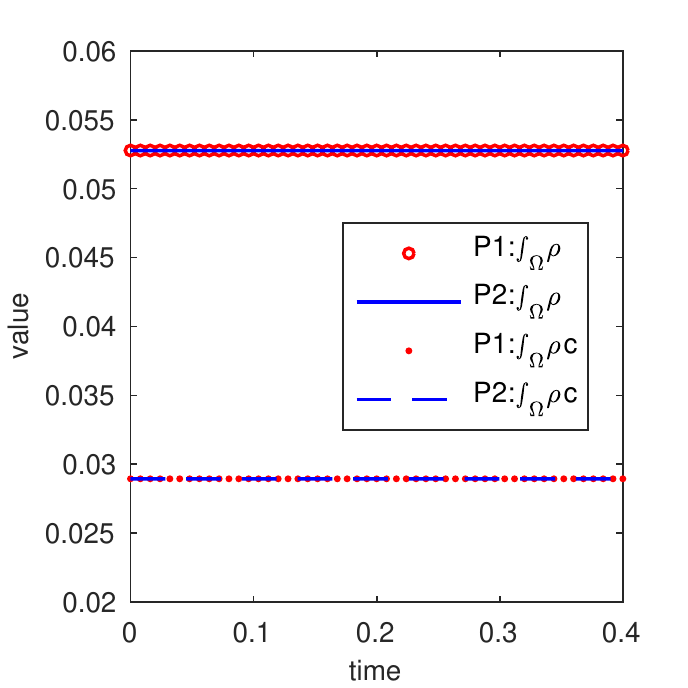}
		\includegraphics[width=3.in]{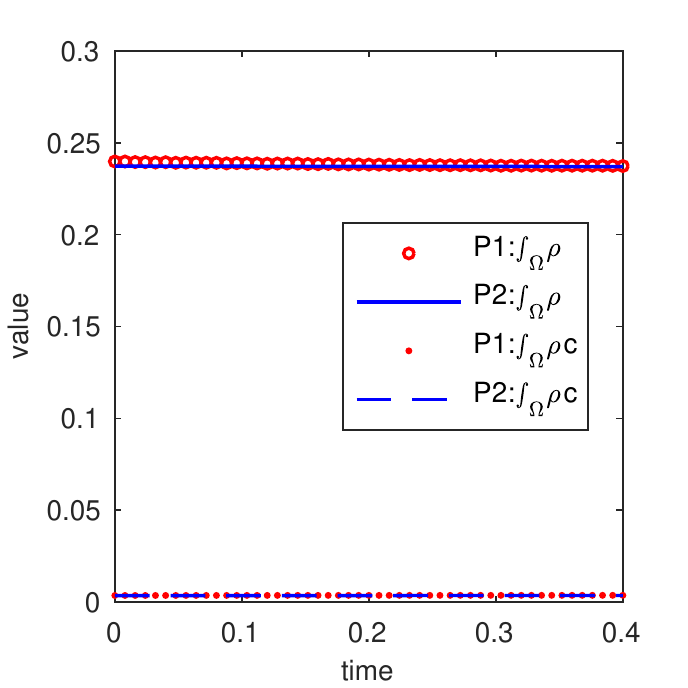}
		\caption{Mass conservation for each component. Left: low density ratio; Right: high density ratio.}
		\label{fig:mass_conserve}
	\end{figure}

	In Fig. \ref{fig:mass_conserve}, we check the total mass convergence of each phase in Lemma \ref{lma:mass_con} for both low and high density ratios. It confirms that P1 and P2  elements could preserve the mass very well in both cases.  
	
\begin{figure}
		\centering
		\includegraphics{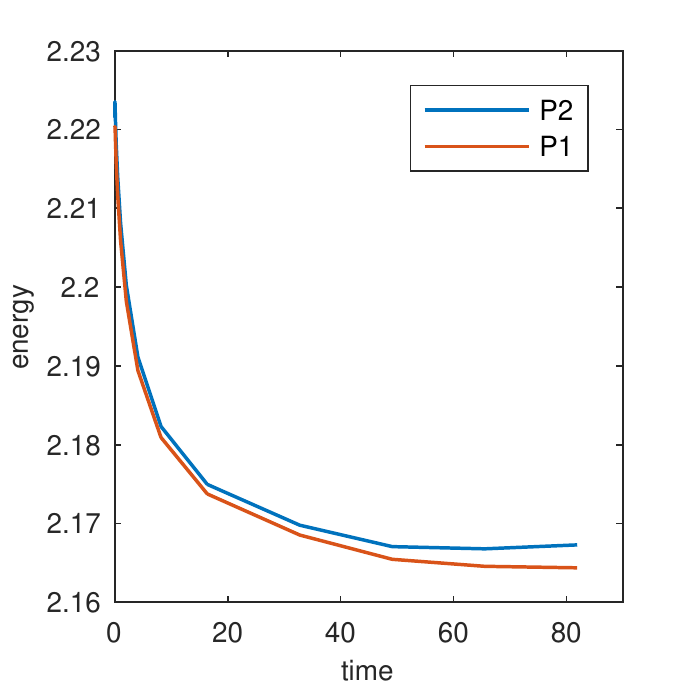}
		\includegraphics{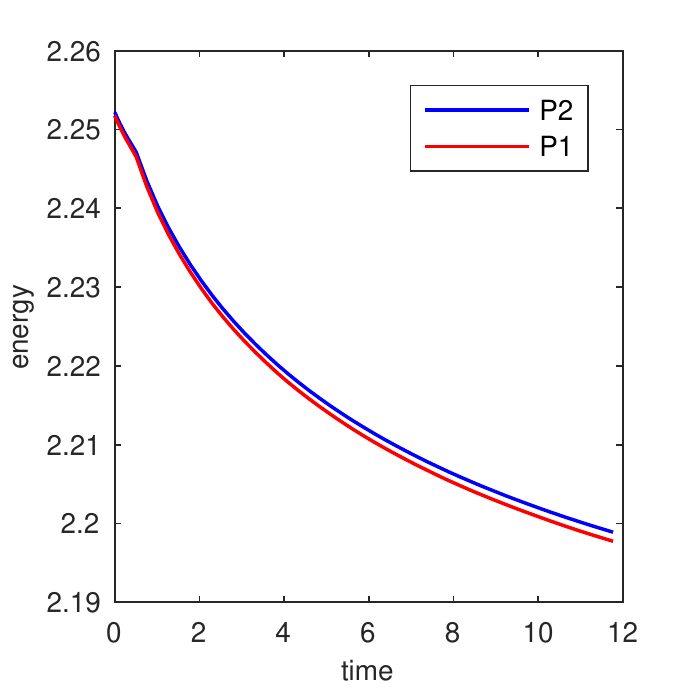}
		\caption{Total energy as a function of time. Left: lower density ratio; Right: high density ratio.}
		\label{fig:energy_decay}
	\end{figure}

	Then we set the wall velocity $\vel_w = (0,0)^T$ to  check the evolution of the total free energy when there is no input energy from outside. It is shown in Fig.\ref{fig:energy_decay} that the free energy decreases over time for both methods and two density ratios, indicating that our schemes are energy stable.

	\subsection{Contact Angle Effect: moving droplet}
	In this example, we show the dynamics of  an oil droplet in water with shear flow. The density ratio is  $\rho_1:\rho_2 = 0.8:1$ and viscosity ratio is $\eta_1:\eta_2 =1:1$.  The other parameters are as follows
	$$Re = 5, \beta = 7.14\times 10^{-3}, M = 2.8\times 10^{-4}, \epsilon = 0.005,\alpha_w = 0.129,$$
	$$ M_{\Gamma} = 5\times 10^8, l_s = 6.667\times 10^{-5}.$$ 
	The domain size is $[0,4]\times[0,0.5]$ with adaptive mesh  and $\Delta t = 4\times 10^{-4}$.	The initial profile is set to be a half circle 
	$$c_0 = 0.5-0.5\tanh\left(\frac{\sqrt{(x-1)^2+y^2}-0.2}{\sqrt{2}\epsilon}\right).$$

\begin{figure}[!ht]
		\centering
		\includegraphics[width=7.in,trim = 50 80 30 80, clip]{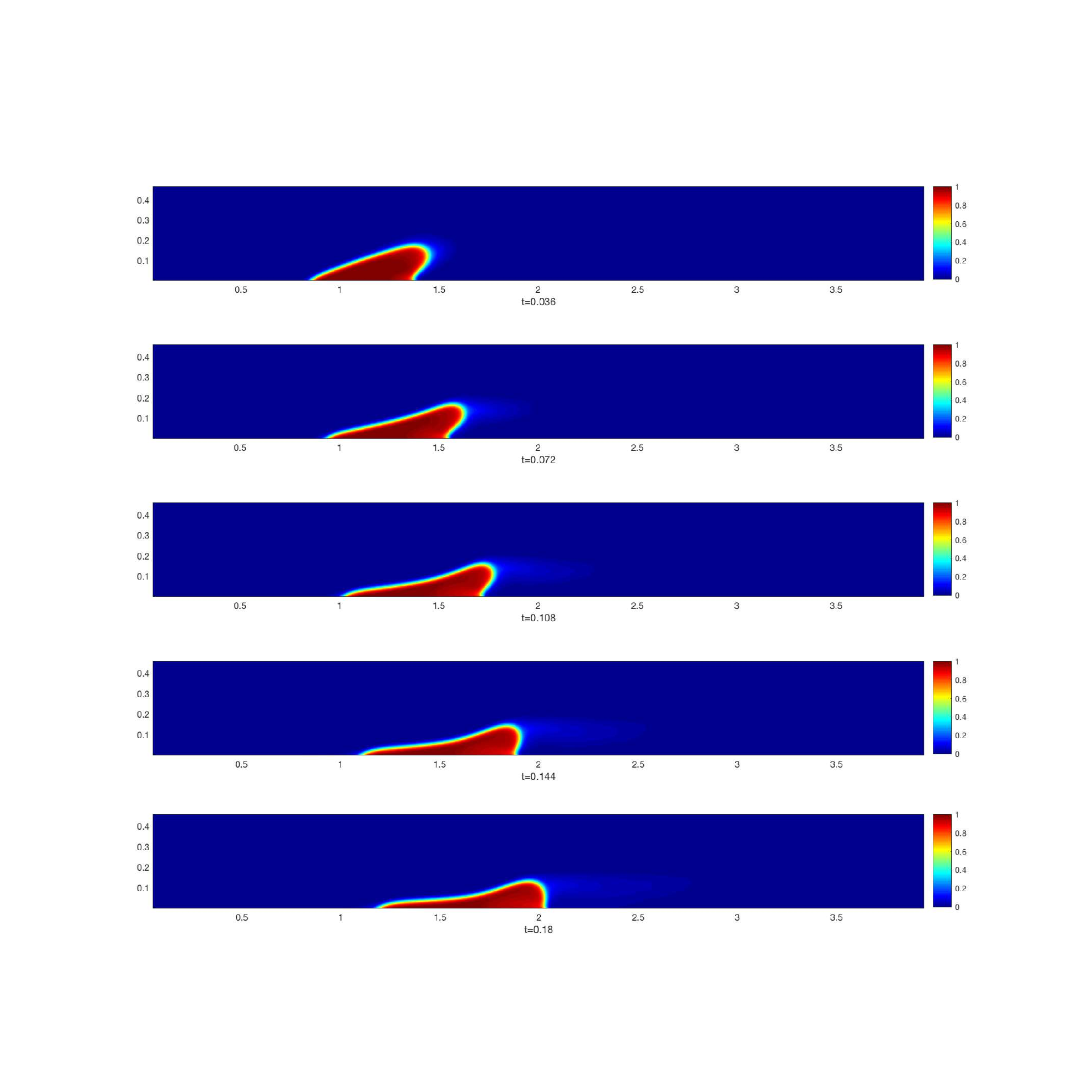}
		\caption{Moving droplet in shear flow with acute static contact angle $\theta_s = 60^{\circ}$.}
		\label{fig:mmovingbubble}
	\end{figure}

\begin{figure}[!ht]
		\centering
		\includegraphics[width=7.in,trim = 50 80 30 80, clip]{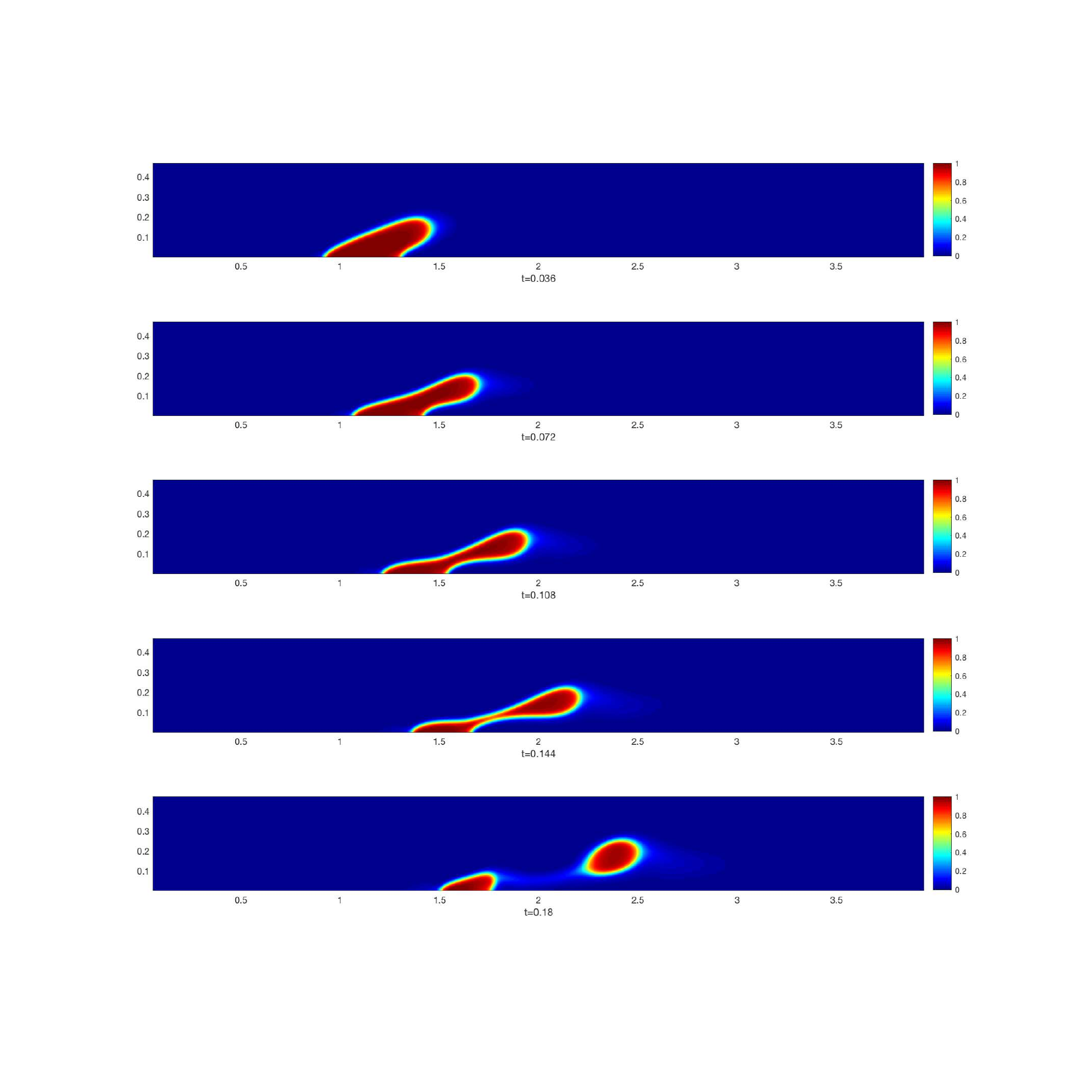}
		\caption{Moving droplet in shear flow with acute static contact angle $\theta_s = 90^{\circ}$.}
		\label{fig:mmovingbubble90}
	\end{figure}

\begin{figure}[!ht]
		\centering
		\includegraphics[width=7.in,trim = 50 80 30 80, clip]{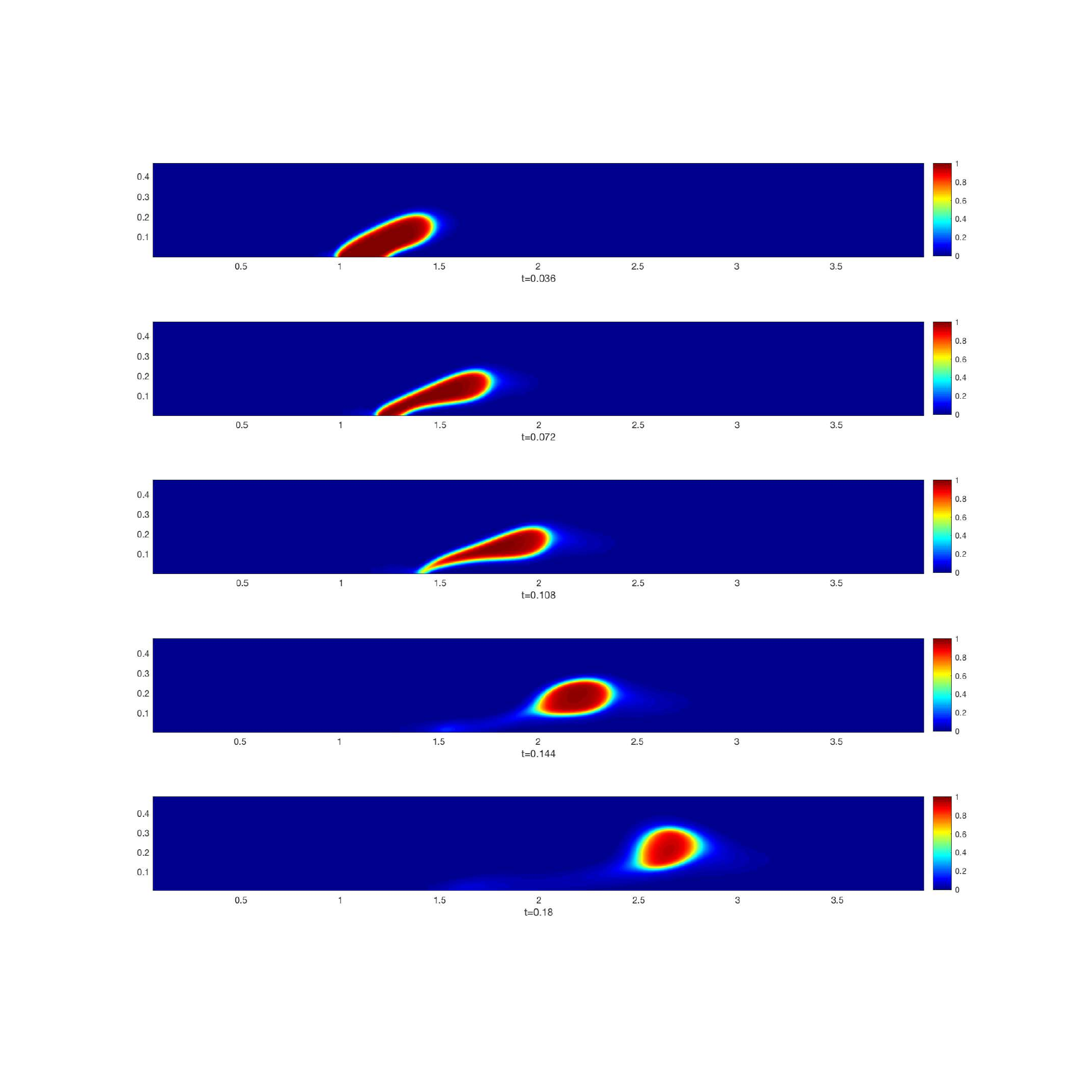}
		\caption{Moving droplet in shear flow with obtuse static contact angle $\theta_s = 120^{\circ}$.}
		\label{fig:mmovingbubble120}
	\end{figure}

\begin{figure}[!ht]
		\centering
		\includegraphics[width=3.in,height=2.5in]{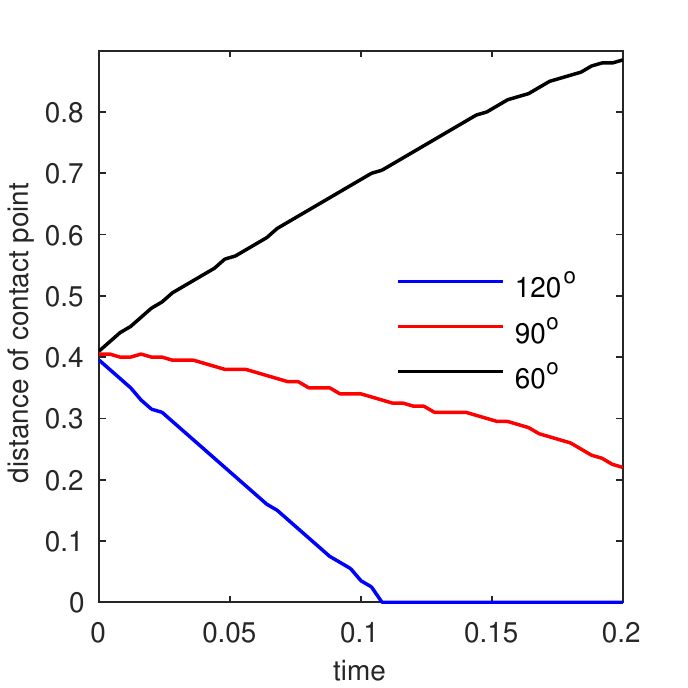}
		\caption{Moving droplet in shear flow: dynamics of the distance between two contact points.}
		\label{fig:mmovingbubblecl}
	\end{figure}

	In Figs. \ref{fig:mmovingbubble} and \ref{fig:mmovingbubble120}, the profiles of droplets under shear flow at different time are presented.  For the acute contact angle cae (Fig. \ref{fig:mmovingbubble} ), the droplet is elongated by  the shear flow force and hydrophilic force on the wall. The distance between two contact points increases over time  (see Fig. \ref{fig:mmovingbubblecl} black curve) as a spreading droplet.  While for the obtuse case (see Fig. \ref{fig:mmovingbubble120}), the hydrophobic force induced the shrink of contact lines on the wall. The distance between two contact points keeps decreasing (see Fig.\ref{fig:mmovingbubblecl} blue curve). With the help of shear force, the droplet eventually detaches from the wall around $t=0.1$ and get stabilized at the center of the flow. 
	When the contact angle is $90^{\circ}$ (Fig. \ref{fig:mmovingbubble90}),
	the competition between wall attraction force and bulk shear force first elongates the droplet and finally breaks the bubble around time  $t= 0.15$. 
	
	\subsection{Large  Density Ratio: Rising Bubble}
	As a last example, we carry out numerical simulation of an air bubble raising in water.
	The density ratio is set to be $\rho_1:\rho_2 = 0.001:1$ and viscosity ratio is $\eta_1:\eta_2 = 0.01:1$. The domain size is $(x,y)\in [0, 0.15]\times[0, 0.15]$ with mesh size $h=1/540$ and timestep $\Delta t = 2\times 10^{-4}$.     
	Parameters are listed as follows
	$$Re = 300, \beta = 0.09, M= 6.67\times 10^{-17}, \epsilon = 0.01,\alpha_w = 100,$$
	$$ M_{\Gamma} = 5\times 10^8, l_s = 0.04.$$ 
	The initial profile is set to be a half circle with radius  0.05  and center at $(0.075,0)$: 
	\begin{equation*}
	c_0=0.5-0.5\tanh\left(\frac{\sqrt{(x-0.075)^2+y^2}-0.05}{\sqrt{2}\epsilon}\right)
	\end{equation*}
	
\begin{figure}
		\centering
		\includegraphics[width=5in,trim = 30 100 30 120, clip]{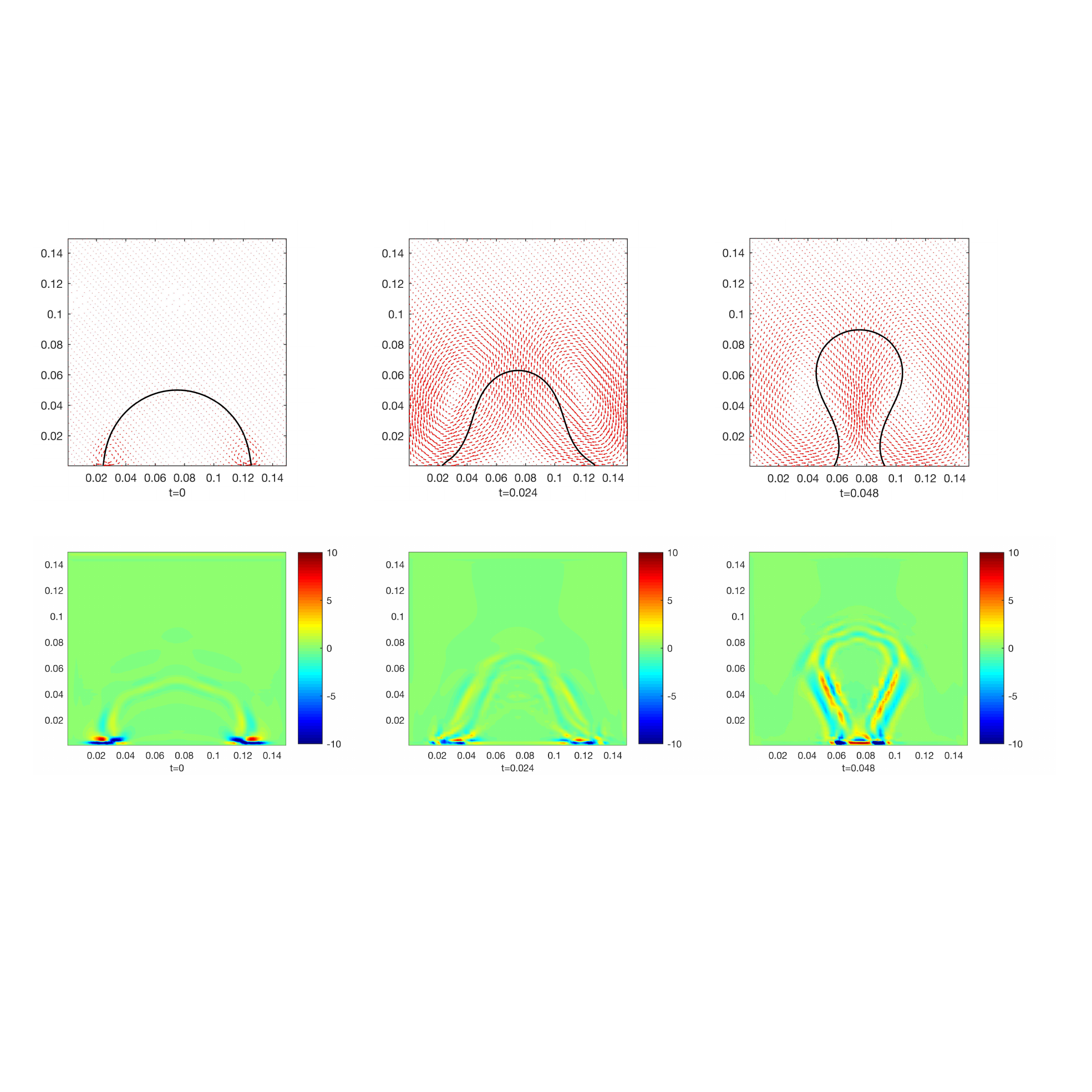}
		\includegraphics[width=5in,trim =30 100 30 120, clip]{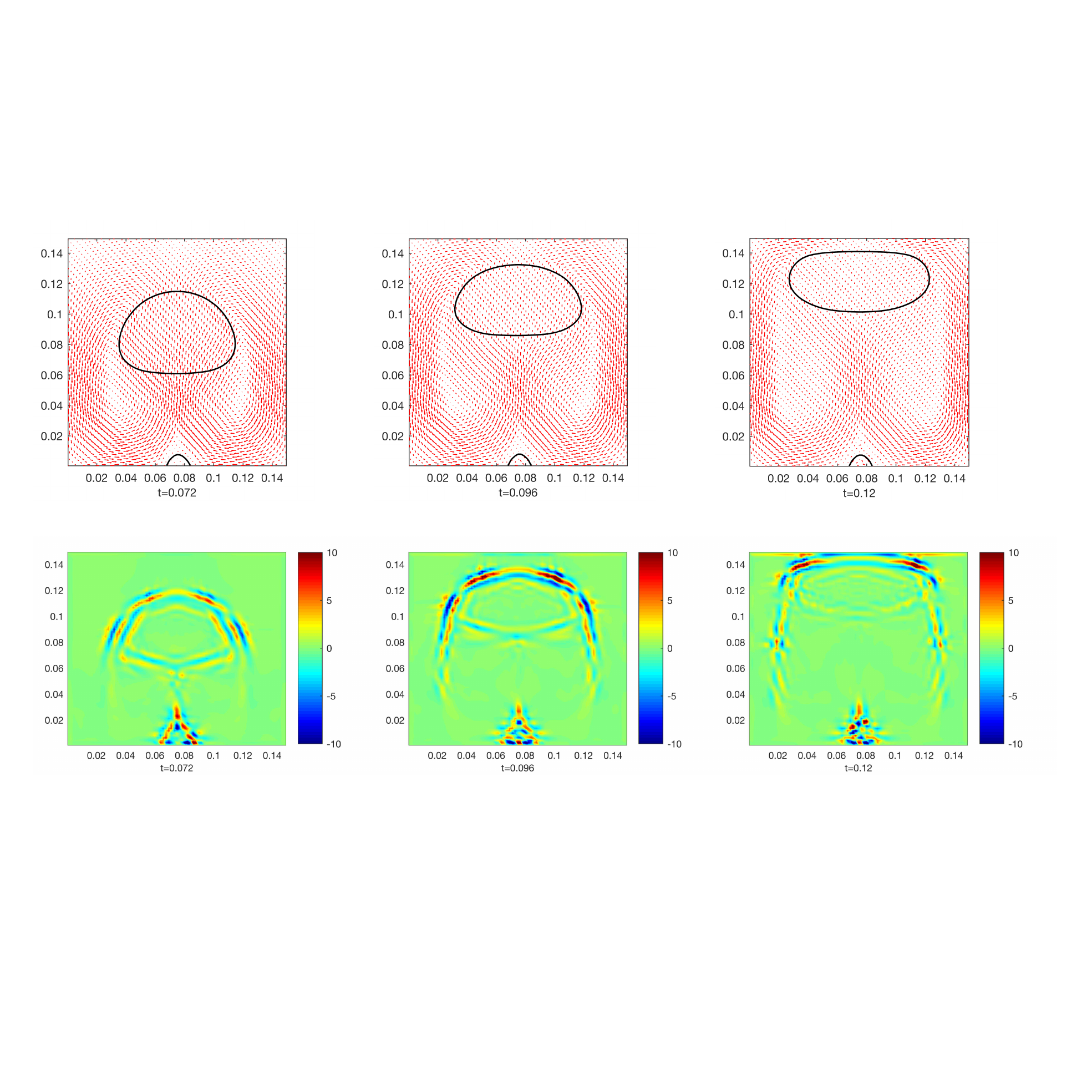}
		\caption{Rising Bubble interface with velocity filed (first and third rows) and $\nabla\cdot \vel$ (second and fourth rows) at different time when $\theta_s = 60^{\circ}$.}
		\label{fig:Rising_Bubble_60}
	\end{figure}

\begin{figure}
		\centering
		\includegraphics[width=5in,trim =30 100 30 120,clip]{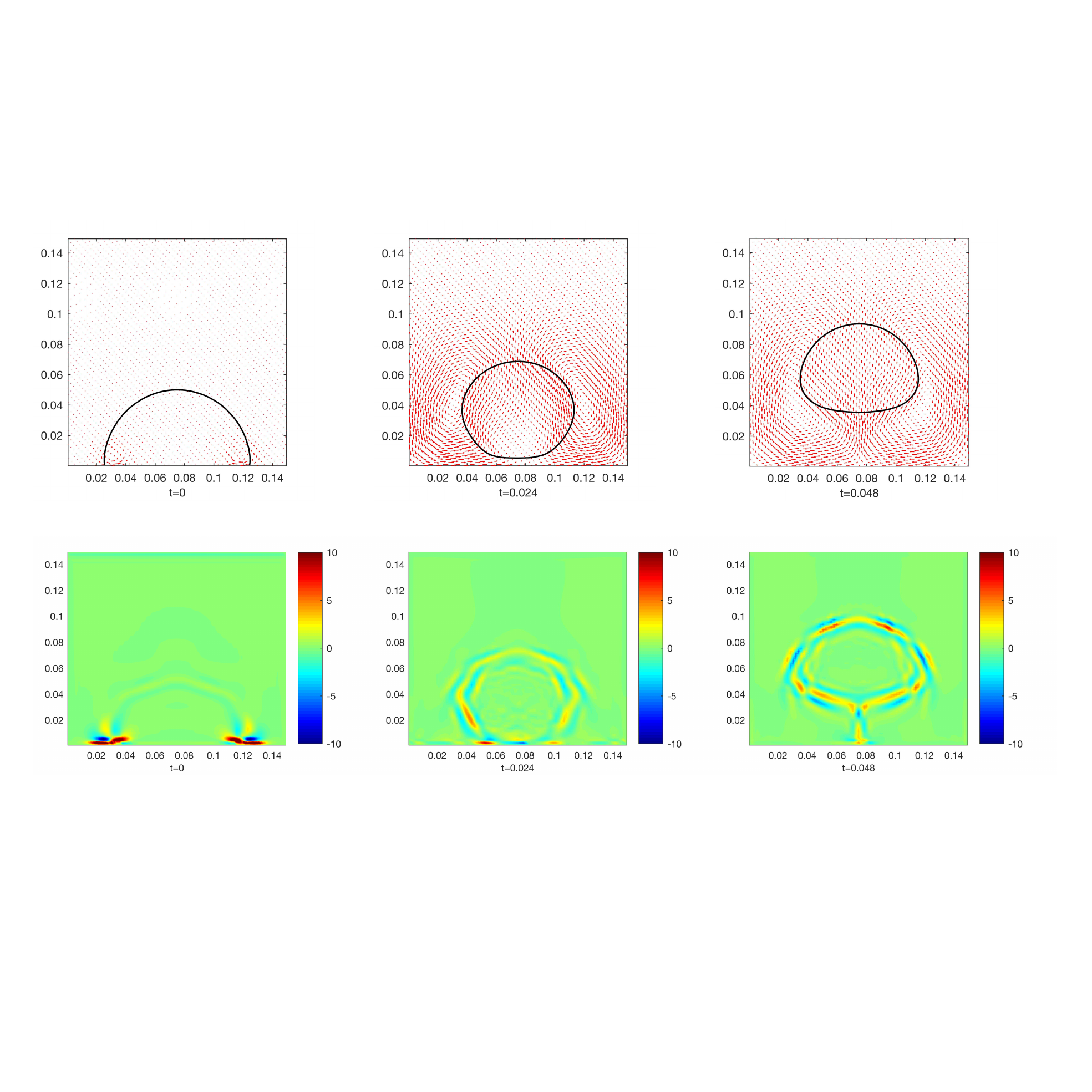}
		\includegraphics[width=5in,trim =30 100 30 120,clip]{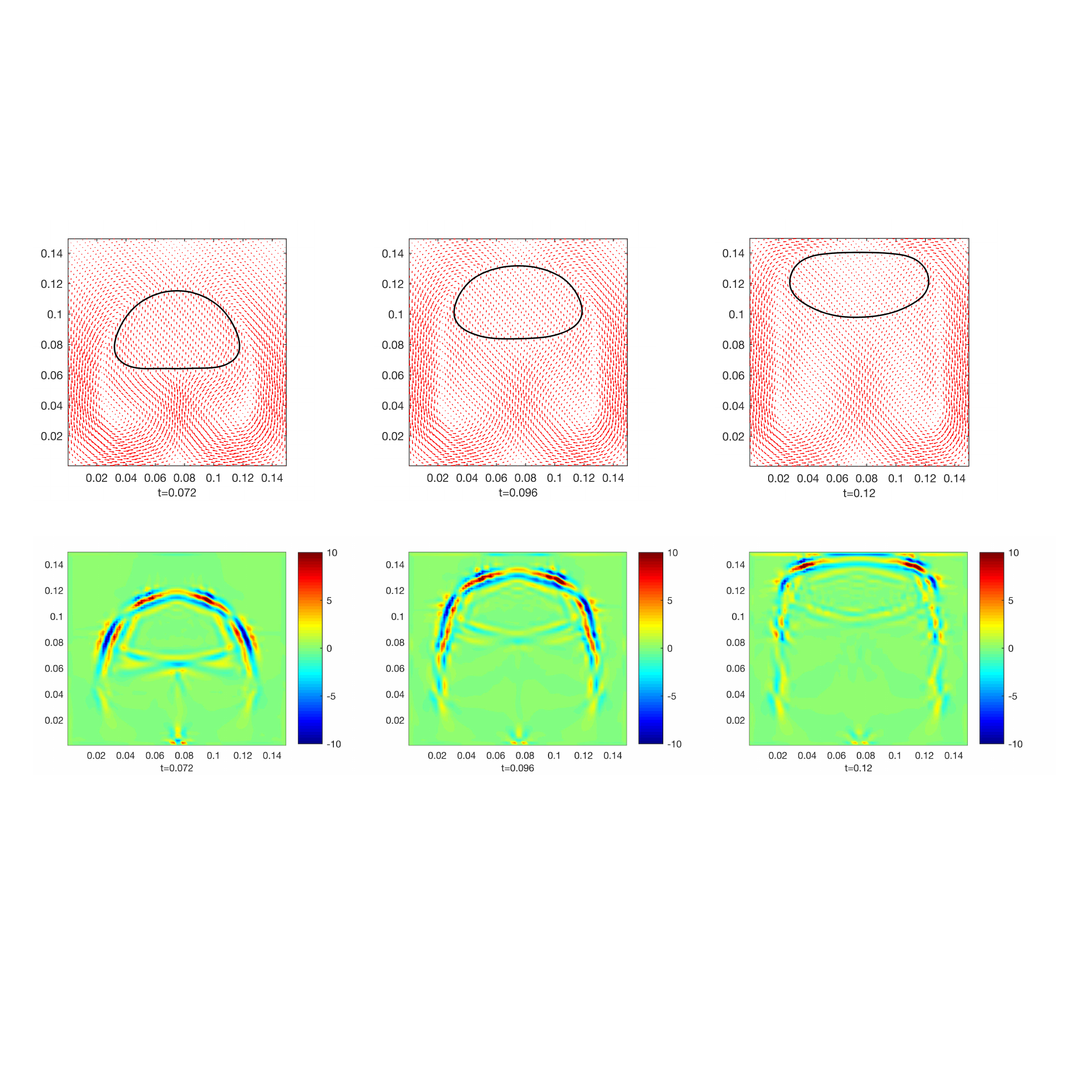}
		\caption{Rising Bubble interface with velocity filed (first and third rows) and $\nabla\cdot \vel$ (second and fourth rows) at different time when $\theta_s = 120^{\circ}$.}
		\label{fig:Rising_Bubble_120}
	\end{figure}

	The snapshots of interfaces with velocity fields and $\nabla\cdot \vel$ profiles for bubbles with acute contact angle $\theta_s = 60^{\circ}$ and obtuse contact angle $\theta_s = 120^{\circ}$ are presented in Figs. \ref{fig:Rising_Bubble_60} and \ref{fig:Rising_Bubble_120}, respectively.  When the angle is acute, the attractive (hydrophilic) force  from the wall competes with the buoyancy force and break the bubble. While for the obtuse case, the wall  repulsive   (hydrophobic) fore enhances the bubble rising under buoyancy force.  The $\nabla\cdot\vel$ profiles confirm that the quasi-impressible property of two-phase fluid with different density only happens around the interface due to the slightly mixing \cite{guo_mass_2017}.  
	
\begin{figure}
		\centering
		\includegraphics[width=3.in,height=2.5in]{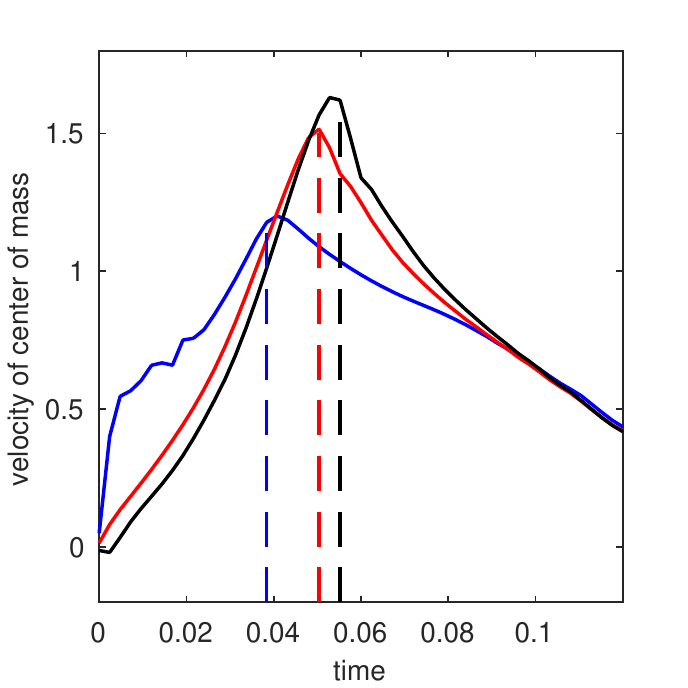}
		\caption{Rising velocity.  Solid lines are rising velocity and vertical dash lines   are time when bubble break or full removed from wall.   Black: $\theta = 60^{\circ}$ with break time $t= 0.0552$;  Red: $\theta = 90^{\circ}$ with break time $t= 0.0504$;  Blue: $\theta = 120^{\circ}$ with removed time $t= 0.0384$.}
		\label{fig:Rising_Bubble_risingvel }
	\end{figure}

	In Fig. \ref{fig:Rising_Bubble_risingvel }, we show the dynamics of rising velocity $V_c =\frac{\int_{\Omega}u_y cd\bx}{\int_{\Omega}  cd\bx}$ of bubble with different static contact angles. The vertical dash lines     are the time when bubbles break ($\theta_s = 60^{\circ}, 90^{\circ}$) or fully detach ($\theta_s = 120^{\circ}$) from wall. It shows that the hydrophobic bubble (Blue line $\theta_s = 120^{\circ}$ ) has a larger acceleration to form a sealing bubble.  At $t= 0.0384$, the bubble fully detaches from the wall. For the hydrophilic bubbles, in the beginning,  the velocity increases slowly due to the competition between the hydrophilic force, the surface tension and the buoyancy force. The bubble  is  stretched into a tear shape  which induces a larger velocity around the narrow neck region (see Fig. \ref{fig:Rising_Bubble_60}).  The maximum velocity is achieved around the break time because the instantaneous response of 
	the surface tension to the large surface deformation.

	\section{Conclusion} \label{section:conclusion}
	In this paper, we first derived the q-NSCH system  for MCLs with variable density by using energy variational method consistently. GNBC for mass-averaged velocity is  obtained during the variation due to the boundary dissipation. 
	
	Then we designed an energy stable $C^0$ finite element scheme to solve the obtained system. We also proved that the fully discrete scheme is   mass conservative for each phase. Thanks to the quasi-incompressible condition with $\Delta p$ term, the finite element space for Navier-Stokes equations do not need to satisfy the Babuska-Brezzi inf-sup condition, as in the case of the pressure stabilization method for the standard Navier-Stokes equations. 
	
	Three examples are investigated numerically. The  Couette flow  test   illustrates the 2nd-order for P1 element and 3rd-order for P2 element  convergence rate in the sense of $L^2$ norm  and the energy decay of the scheme.  The contact angle effect on the droplet is illustrated by moving droplet in shear flow. Finally, a rising bubble is simulated to confirm the ability of our scheme to handle large density ratio and the quasi-incompresiblity only happens around the interface.


	\noindent\section*{Acknowledgments}
	
	This work was partially supported by the NSFC (grant numbers 11771040, 11861131004), NSERC (CA) and the Fields Institute. Lingyue Shen was partially supported by the Chinese Scholarship Council for studying at the University of Dundee.

	
	
	\bibliographystyle{elsarticle-num} 
	\bibliography{Mybib}

	\newpage
	\appendix
	\section{Energy Variation Details}\label{energyvariation}
	
	For the first term $I_1$ in (\ref{dEdt}), using the last two equations in Eq.\eqref{assumption} yields
	\be\label{I_1a}
	I_1 &=&\frac{d}{dt}\int_{\Omega}\frac{\rho|\vel|^2}{2}d\bx \nonumber\\
	&=&\int_{\Omega}\frac 1 2 \frac{\partial\rho}{\partial t}|\vel|^2d\bx +\int_{\Omega} \rho\frac{\partial \vel}{\partial t}\cdot\vel d\bx\nonumber\\
	&=&\int_{\Omega}\frac 1 2 \frac{\partial\rho}{\partial t}|\vel|^2d\bx +\int_{\Omega} \rho\frac{D \vel}{D t}\cdot\vel d\bx-\int_{\Omega} \left(\rho \vel\cdot\nabla \vel\right)\cdot \vel d\bx\nonumber\\
	&=&\int_{\Omega}\frac 1 2 \frac{\partial\rho}{\partial t}|\vel|^2d\bx +\int_{\Omega} \rho\frac{D \vel}{Dt}\cdot\vel d\bx+\int_{\Omega} \nabla\cdot(\rho\vel)\frac{|\vel|^2}2d\bx\nonumber\\
	&=&\int_{\Omega} \left(\nabla \cdot (\bsigma_\eta + \bsigma_ c) \right) \cdot\vel d\bx +\int_{\Omega}  p \left( \frac{1}{\rho^2} \frac{d  \rho}{d c}\nabla\cdot \bj_{ c}-\nabla\cdot \vel  \right)d\bx \nonumber\\
	&=&\int_{\Omega} \left(\nabla \cdot (\bsigma_\eta + \bsigma_ c) \right) \cdot\vel d\bx +\int_{\Omega}  p \left( -\alpha\nabla\cdot \bj_{ c}-\nabla\cdot \vel  \right)d\bx \nonumber\\
	&=&-\int_{\Omega}(\bsigma_{\eta}:\nabla\vel +\bsigma_{ c}:\nabla\vel)d\bx+\int_{\Omega}\nabla (\alpha  p)  \cdot \bj_{ c}d\bx-\int_{\Omega} p  \nabla\cdot \vel d\bx\nonumber\\
	&&+\int_{\partial\Omega_w}((\bsigma_{\eta}+\bsigma_{ c})\cdot\bn)\cdot\vel_{\tau} dS.
	\ee
	where we have introduced a Lagrangian multiplier $p$ with respect to the constraint (\ref{quasi_incom3}) and have used the boundary condition $\vel\cdot\bn = 0$ and $\bj_{ c}\cdot\bn=0$.

	For the second term $I_2$ in (\ref{dEdt}), using the first equation in Eq.\eqref{assumption} and last two boundary conditions in Eq. \eqref{assumption_bd} yields
	\be\label{I_2a}
	I_2& =&\frac{d}{d t}\int_{\Omega}\rho\lambda_{ c}\left(G( c)+\frac{\gamma^2}{2}|\nabla c|^2\right)d\bx\nonumber\\
	&=&\int_{\Omega}\rho\lambda_{ c}\frac D {Dt}\left(G+ \frac {\gamma^2}2|\nabla c|^2 \right)d\bx\nonumber\\
	&=&\int_{\Omega}\rho\lambda_{ c}\frac {dG} {d c}\frac{D c}{Dt}d\bx+ \int_{\Omega} \rho\lambda_{ c} \gamma^2\left(\nabla c\cdot\frac{D}{Dt}(\nabla c) \right)d\bx\nonumber\\
	&=&\int_{\Omega}\rho\lambda_{ c}\frac {dG} {d c}\frac{D c}{Dt}d\bx+ \int_{\Omega} \rho\lambda_{ c} \gamma^2\left(\nabla c\cdot\left(\frac{\partial}{\partial t}(\nabla c)+(\vel\cdot\nabla)(\nabla c)\right) \right)d\bx\nonumber\\
	&=&\int_{\Omega}\rho\lambda_{ c}\frac {dG} {d c}\frac{Dc}{Dt}d\bx+ \int_{\Omega} \rho\lambda_{ c} \gamma^2\left(\nabla c\cdot\left(\nabla\frac{\partial c}{\partial t}\right)\right)d\bx +\int_{\Omega}\rho\lambda_{ c} \gamma^2(\partial_i c u_j\partial^2_{ji} c)d\bx\nonumber\\
	&=&\int_{\Omega}\rho\lambda_{ c}\frac {dG} {d c}\frac{D c}{Dt}d\bx+ \int_{\Omega} \rho\lambda_{ c} \gamma^2\left(\nabla c\cdot\left(\nabla\frac{\partial c}{\partial t}\right)\right)d\bx \nonumber\\
	&&+\int_{\Omega}\rho\lambda_{ c} \gamma^2(\partial_i c \partial_i(u_j\partial_{j} c) -\partial_i c \partial_j c\partial_iu_j )d\bx\nonumber\\
	&=&\int_{\Omega}\rho\lambda_{ c}\frac {dG} {d c}\frac{Dc}{Dt}d\bx+ \int_{\Omega} \rho\lambda_{ c} \gamma^2\nabla c\cdot\nabla\left(\frac{D c}{D t}\right)d\bx-\int_{\Omega}\rho\lambda_{ c} \gamma^2(\nabla c \otimes\nabla c):\nabla\vel d\bx\nonumber\\
	&=&\int_{\Omega}\rho\lambda_{ c}\frac {dG} {d c}\frac{D c}{Dt}d\bx- \int_{\Omega}\nabla\cdot \left(\rho\lambda_{ c} \gamma^2\nabla c\right)\left(\frac{D c}{D t}\right)d\bx-\int_{\Omega}\rho\lambda_{ c} \gamma^2(\nabla c \otimes\nabla c):\nabla\vel d\bx\nonumber\\
	&& +\int_{\partial\Omega_w}\rho \lambda_{ c}\gamma^2\partial_n c \frac{D_{\Gamma} c}{Dt}dS\nonumber\\
	&=&\int_{\Omega}\rho\mu \frac{D c}{D t} d\bx-\int_{\Omega}\lambda_{ c} \gamma^2(\rho\nabla c \otimes\nabla c):\nabla\vel d\bx +\int_{\partial\Omega_w}\rho \lambda_{ c}\gamma^2\partial_n c \frac{D_{\Gamma} c}{Dt}dS\nonumber\\
	&=&-\int_{\Omega}\mu\nabla\cdot\bj_{ c} d\bx-\int_{\Omega}\lambda_{ c} \gamma^2(\rho\nabla c \otimes\nabla c):\nabla\vel d\bx+\int_{\partial\Omega_w}\rho \lambda_{ c}\gamma^2\partial_n c \frac{D_{\Gamma} c}{Dt}dS\nonumber\\
	&=&\int_{\Omega}\nabla\mu\cdot\bj_{ c} d\bx-\int_{\Omega}\lambda_{ c} \gamma^2(\rho\nabla c \otimes\nabla c):\nabla\vel d\bx+\int_{\partial\Omega_w}\rho \lambda_{ c}\gamma^2\partial_n c \frac{D_{\Gamma} c}{Dt}dS\nonumber\\
	\ee
	where $\mu = \lambda_{ c} \left(\frac{dG}{d c}-\frac{1}{\rho}\gamma^2\nabla\cdot(\rho\nabla c)\right)$.
	\section{Proof of Lemma \ref{transportlma}}
	\noindent\textbf{Proof}: \begin{eqnarray}
	&&\frac{d}{dt}\int_{\Omega}\left(\rho(\bx,t) f(\bx,t)\right)d\bx\nonumber\\
	&=& \int_{\Omega}\frac{\partial\rho}{\partial t}f +\rho\frac{\partial f}{\partial t}d\bx\nonumber\\
	&=&-\int_{\Omega}\nabla\cdot(\vel\rho)fd\bx +\int_{\Omega} \rho \frac{\partial f}{\partial t}d\bx \nonumber\\
	&=& \int_{\Omega}\rho (\frac{\partial f}{\partial t}+\vel\cdot \nabla f)d\bx\nonumber\\
	&=& \int_{\Omega}\rho\frac{Df}{Dt}d\bx.
	\end{eqnarray}

	\section{Tensor Calculation}\label{tensor_calculation}
	\begin{eqnarray}
	&& (\nabla\vel +(\nabla\vel)^T):\nabla\vel -\frac{2}{3}(\nabla\cdot\vel)^2\nonumber\\
	&=&\sum_{i,j=1,2,3}(\partial_iu_j +\partial_ju_i)\partial_ju_i -\frac 2 3  \sum_{i=1,2,3} (\partial_i u_i)^2 \nonumber\\
	&=& \sum_{i=1,2,3}2(\partial_iu_i)^2+ \sum_{i<j}2(\partial_iu_j\partial_ju_i) +\sum_{i\neq j} (\partial_iu_j)^2-\frac{2}{3}\sum_{i=1,2,3} (\partial_i u_i)^2 - \frac{2}{3}\sum_{i<j}2(\partial_iu_i\partial_ju_j)\nonumber\\
	&=& \sum_{i<j} (\partial_iu_j+\partial_ju_i)^2+\frac{4}{3}(\sum_{i=1,2,3}(\partial_iu_i)^2-\sum_{i<j} (\partial_iu_i\partial_ju_j))\nonumber\\
	&=&  \sum_{i<j} (\partial_iu_j+\partial_ju_i)^2+\frac{2}{3}  \sum_{i<j} (\partial_iu_i-\partial_ju_j)^2
	\end{eqnarray}
	
	
	
	
	
\end{document}